\documentclass[12pt]{iopart}
\usepackage{amsopn,times,graphicx,amsthm}
\usepackage[percent]{overpic}
\usepackage{varioref}
\usepackage{placeins}
\usepackage[matrix,arrow]{xy}
\usepackage{array}
\usepackage{longtable}
\usepackage{dcolumn}
\usepackage{booktabs}
\usepackage{epic,eepic}
\usepackage{fp}
\usepackage{units}
\usepackage{url}
\def\marginpar#1{}   

\let\lbl=\label
\def\label#1{\lbl{#1}\ifinner\else\marginpar{\ref{#1} #1}\ignorespaces\fi}

\newcommand\Rop{\operatorname{Rop}}

\newcommand\PRop{\operatorname{Rop}_p}

\renewcommand{\phi}{\varphi}

\theoremstyle{definition}

\newcommand{\FileCount}{943}

\newcommand{\CompositeKnotCount}{544}
\newcommand{\CompositeKinkCount}{542}
\newcommand{\CompositeStraightCount}{534}

{\makeatletter
 \gdef\xxxmark{%
   \expandafter\ifx\csname @mpargs\endcsname\relax 
     \expandafter\ifx\csname @captype\endcsname\relax 
       \marginpar{xxx}
     \else
       xxx 
     \fi
   \else
     xxx 
   \fi}
 \gdef\xxx{\@ifnextchar[\xxx@lab\xxx@nolab}
 \long\gdef\xxx@lab[#1]#2{{\bf [\xxxmark #2 ---{\sc #1}]}}
 \long\gdef\xxx@nolab#1{{\bf [\xxxmark #1]}}
}

\newcommand{\Wr}{\operatorname{Wr}}
\newcommand{\Cr}{\operatorname{Cr}}

\begin{document}

\title{Shapes of Tight Composite Knots}
\date{November 3, 2010; Revised: \today}

\author{Jason Cantarella$^1$, Al LaPointe$^1$ and Eric J.~Rawdon$^2$}
\address{$^1$ Department of Mathematics, University of Georgia, Athens, GA 30602, USA}
\address{$^2$ Department of Mathematics, University of St.~Thomas, Saint Paul, MN 55105, USA}

\begin{abstract}
  We present new computations of tight shapes obtained using the
  constrained gradient descent code \texttt{ridgerunner} for
  \CompositeKnotCount\ composite knots with 12 and fewer crossings,
  expanding our dataset to \FileCount\ knots and links. We use
  the new data set to analyze two outstanding conjectures about tight
  knots, namely that the ropelengths of composite knots are at least
  $4\pi-4$ less than the sums of the prime factors and that the writhes
  of composite knots are the sums of the writhes of the prime factors.
\end{abstract}   


\maketitle

\section{Introduction}

When people are given a piece of rope, they almost instinctively tie
it into a knot and pull it tight.  But what, exactly, is the structure
of that tight knot?  The last decade has seen great progress in
analyzing tight configurations.  Researchers have focused on
mathematical knots, that is closed loops forming different topological
knot types, with the rope modeled as a non self-intersecting tube about a
smooth (usually $C^2$ or $C^{1,1}$) space curve. One can then define
the
\textit{ropelength} of the curve to be the quotient of its length and the
maximal radius (or \textit{thickness}) of a non self-intersecting tube
about the curve.  Alternately, the ropelength is the minimal
centerline length of a unit radius tube without self-intersections
forming the given knot type.

Configurations that minimize the ropelength within a given knot type
are called \textit{tight} or \textit{ideal}.  These configurations
have been used to predict the relative speed of DNA knots under gel
electrophoresis~\cite{MR99b:92032}, the pitch of double helical
DNA~\cite{mbms}, the average values of different spatial measurements
of random knots~\cite{MR1981020}, and the breaking points of
knots~\cite{1367-2630-3-1-310}. They also provide a model for the
structure of a class of subatomic particles known as
glueballs~\cite{kepglue}.  Another way to think about the tight knot
problem is to see it as a packing problem akin to Kepler's
Conjecture. In this version, instead of packing individual spheres
into a volume, we are packing an entangled tube into a small volume.

Finding analytical solutions to the tight knot problem is difficult;
we know the tight configurations for only some specialized classes of
links~\cite{MR2003h:58014}.  For even the simplest non-trivial knot,
the trefoil $3_1$, there is no analytic solution for the tight
configuration.  Instead, researchers rely on computer simulations to
approximate tight configurations by polygons minimizing an
appropriately discretized version of
ropelength~\cite{rawdonthesis,MR1749504}.  Such calculations yield
upper bounds for the minimal ropelength of knots and links. Combining
computer simulations with theoretical work, we know that the minimal
ropelength of the trefoil is between $31.32$~\cite{MR2207788} and
$32.74317$~\cite{MR2471403}.

The quality of computer approximations of tight knots has increased
immensely over time.  Originally, simple techniques such as simulated
annealing were used to determine approximately tight
configurations~\cite{MR99b:92032,MR2034393}.  Later, P\'ieranski wrote
the \texttt{SONO} (Shrink On No Overlaps) software
\cite{MR1702020,tightopen,baranska,gordian}, implementing a
gradient-like algorithm that shortens the length of the polygon and
then pushes pairs of vertices apart when they create
self-intersections in a tube about the polygon.  Maddocks's group has
used a different approach, implementing simulated annealing on biarc
curves to minimize the ropelength of smooth curves directly
\cite{smutnyphd,maddocksbook}.

Over the past few years, our group has developed the knot-tightening
code \texttt{ridgerunner}, which implements a constrained
gradient-descent algorithm that minimizes the length of a polygon
subject to a family of constraints which define an embedded tube
around the polygon~\cite{acprtightening}. The algorithm projects the
gradient of the length of the polygon onto the subspace of motions
that preserve the integrity of the unit-radius embedded tube
about the polygon; we prove
in~\cite{acprtightening} that a polygon is ropelength-critical when
this projection vanishes. We use this fact to define a quality measure
for approximately tight configurations: the \textit{residual} of such
a configuration is the fraction of the ($L^2$) norm of the gradient of
length after projection onto the constraint.

Since tight knot configurations are so useful in the sciences, it is
desirable to have a complete catalogue of tight knot shapes. The first
step is to assemble a table of knot types. Knot tabulation has a long
history, stretching back to the 19th century knot tables of Tait and
Kirkman~\cite{tait,kirkman} and very large tables of knots have been generated by
computer~\cite{MR2105106,MR1646740}. Using such tables, we have computed approximately
tight shapes for 379 knots and links with 10 and fewer
crossings~\cite{acprtightening}. However, these tables are incomplete
in a certain sense: they only contain ``prime'' knots. To understand
primeness for knots, consider the knot ``product'' defined by splicing
two knots $K_1$ and $K_2$ together. This is called the \emph{connect
  sum} of the knots and denoted $K_1 \# K_2$. We say that $K$ is prime
if $K = K_1 \# K_2$ only when $K_1$ or $K_2$ is an unknotted loop in
analogy to the idea of primality for natural numbers, where $n$ is
prime if and only $n = km$ implies $k$ or $m$ is equal to $1$.  The
standard knot tables, and our work in~\cite{acprtightening} include
only prime knots.

But while prime knots and links are mathematically convenient, there
is no reason to expect that the knots and links which occur in
scientific applications will be prime. For this reason, we have
continued our work to tabulate and tighten composite knots. In this
paper, we present the results of a large-scale computation of
approximately tight shapes for composite knots with 12 and fewer
crossings, covering \CompositeKnotCount\ knot types. We report the
ropelengths of these shapes in Tables~\ref{ByKnotTable
  0}-\ref{ByKnotTable 5} on pages~\pageref{ByKnotTable
  0}--\pageref{ByKnotTable 5}. We computed these shapes with
\texttt{ridgerunner}, generating starting configurations by splicing
together the approximately tight configurations of prime knots
from~\cite{acprtightening}. It is an open problem whether the crossing
number of a composite knot is the sum of the crossing numbers of its
prime factors; if this is true, our list of composite knots covers all
the composites with 12 and fewer crossings. The quality of these
computations is measured by their resolution (the number of vertices
per unit ropelength) and their residual (the fraction of the
tightening force which is not balanced by contact forces on the
tube). All of our knots have resolution at least 8, and almost
all of our knots have residuals
between 0.01 and 0.001. The residual for each shape is reported in
Tables~\ref{ResidualTable 0}--\ref{ResidualTable 5} on
pages~\pageref{ResidualTable 0}--\pageref{ResidualTable 5}.  Together
with our earlier work, this brings the total number of computed
approximately tight shapes to \FileCount\ knot and link types.

Katritch, Olsen, Pieranski, Dubochet and Stasiak made the first
computations of tight composite knots, reported in \emph{Nature} in
1997~\cite{kops}. In that paper, they observed several interesting
phenomena. First, they noticed that the 3D average writhing number of
tight configurations appeared to be additive under connect sum:
$\Wr(K_1) + \Wr(K_2) = \Wr(K_1 \# K_2)$. This was a particularly
striking observation since the writhe is a shape invariant and
\emph{not} a topological invariant and there is no reason to believe
that the tight shapes of $K_1$ and $K_2$ would be exactly repeated in
the tight configuration of $K_1 \# K_2$. We find that for the vast
majority of knots, this conjectured equation holds to a remarkable
degree of accuracy. However, we have found a small number of anomalous
composites where the conjecture seems to fail.

Another phenomenon noted in~\cite{kops} was that the ropelength of the
tight configuration of $K_1 \# K_2$ was shorter than the sum of the
minimum ropelengths of $K_1$ and $K_2$: intuitively, one could save a certain
amount of rope by splicing. They conjectured that the amount of rope
saved was at least $4\pi-4$ (this is true for the links
in~\cite{MR2003h:58014}). Our computations support this conjecture,
although we find many cases where the amount of rope saved is very close
to the conjectured $4\pi - 4$.

\section{Methodology}

\subsection{Tabulation of composite knots} 
It is a classical theorem of Schubert~\cite{MR0031733} that every
composite knot has a unique decomposition into an unordered list of
prime summands. From this perspective, it would seem that tabulating
composite knots must be easy: one should take a knot table and form
all subsets of the table (allowing repeats). Unfortunately, the
situation is not quite this simple.  The traditional tables of prime
knots are computed only up to symmetry. For each entry in the tables,
there are actually one, two, or four distinct knot types associated
with $K$, depending on the symmetry properties of $K$.

It is easiest to see this with respect to chirality, a topic of great
familiarity in the sciences.  The trefoil knot $3_1$, for example, is
\textit{chiral} in that one cannot deform a ``right-handed'' trefoil
to its mirror image. Thus the mirror image of the trefoil is actually
a member of a different knot type, denoted $3_1^m$ and called the
``left-handed'' trefoil. On the other hand, the figure-8 knot $4_1$
can be deformed to its mirror image, so it is called
\textit{amphichiral}.  Most knots are chiral so they are not
topologically equivalent to their mirror images, but some are
amphichiral (for example, 20 of the 249 prime knots with 10 or fewer
crossings are amphichiral).

A \textit{reversible} knot type is one that can be deformed to itself
but with the opposite orientation along the curve. A symmetric
configuration of the trefoil, for instance, can be reversed by
rotating it by 180 degrees. Only 36 of the 249 prime knots with 10 or
fewer crossings are non-reversible (note that the ratio of
non-reversible to reversible knot types increases with crossing
number), the simplest of which is the $8_{17}$ knot. We denote the
reverse of a knot type $K$ by $K^r$, so the reverse of the $8_{17}$
knot is $8_{17}^r$. From a physical standpoint, the reversibility of a
knot type could be as important as its chirality, for instance when
the knot represents a flux tube or strand of DNA and thus has a
natural orientation.

As a result, there are four obvious classes of knots:
chiral/non-reversible (no symmetry), chiral/reversible (invertible symmetry,
although Conway~\cite{MR0258014} calls this \emph{reversible}
symmetry), amphichiral/non-reversible ((+) amphichiral symmetry, although
Conway~\cite{MR0258014} calls this \emph{invertible} symmetry), and amphichiral/reversible (full symmetry). In
addition, there is a class of \textit{negative amphichiral} knot types
which are not equivalent to either their reverses or their mirror
images, but \textit{are} equivalent to the reverse of their mirror
images. This symmetry type does not fit neatly into the classification
above. Luckily, these knots are rare among knots of low crossing
number. A summary of the five different knot symmetries are given in
Table~\ref{tab:symmetrytypes}.

\begin{table}[ht]
	\begin{center}
  \begin{tabular}{ccccc}
    \toprule
    Class & Amphichiral & Reversible & Isotopy Types & Example(s)\\
    \midrule
    No symmetry & No & No & 4 & $9_{32}, 9_{33}$ \\
    (-) amphichiral symmetry & - & - & 2 & $12_{427}$ \\
    invertible
    symmetry & No & Yes & 2 & $3_1$ \\
    (+) amphichiral
    symmetry  &  Yes & No & 2 & $8_{17}$ \\
    full symmetry & Yes & Yes & 1 & $4_1$ \\
    \bottomrule
  \end{tabular}
  \end{center}
  \medskip
  \caption{The five standard symmetry types for a knot type.}
  \label{tab:symmetrytypes}
\end{table}

For tightening computations on prime knots, these fine distinctions
are usually immaterial. For instance, although the trefoil knot
($3_1$) is not isotopic to the mirror trefoil ($3_1^m$) we know that
any tight configuration of $3_1$ is a rigid reflection of a tight
configuration of $3_1^m$. Hence both of these knot types have the same
minimum ropelength. We note that other geometric invariants which are
sensitive to chirality, such as the \emph{average writhing number},
will be different for the tight configuration of each knot type.

However, when considering composites, symmetries make a real
difference in the shapes of tight knots. The connect sum of two
trefoils $3_1 \# 3_1$ (called the \emph{granny knot}) is not only a
different knot type from the connect sum $3_1 \# 3_1^m$ (the
\emph{square knot}) but it also has a different minimum ropelength. On
the other hand, the \emph{mirror granny} knot $3_1^m \# 3_1^m$ is a
different knot type than the granny knot $3_1 \# 3_1$ but has the same
minimum ropelength, while the knot \emph{mirror square} knot $3_1^m \#
3_1$ has the same knot type, and thus the same minimum ropelength, as
the square knot $3_1 \# 3_1^m$. It turns out there are various
possibilities when one takes the connect sum of two knots depending on
their symmetry types, with a further simplification when the summands
are related to one another by a symmetry (such as in the case
above). These possibilities are summarized in
Table~\ref{tab:ktcounts}.

\begin{table}[ht]
\begin{center}
\begin{tabular}{c|ccccc}
\toprule
$\#$ & None & (-) Amphichiral & Reversible & (+) Amphichiral & Full \\
\midrule
None & 16 (4) & & & & \\
(-) Amphichiral & 12 (2) & 9 (2) & & & \\
Reversible & 8 (2) & 6 (1) & 4 (2) & & \\
(+) Amphichiral & 8 (2) & 6 (1) & 4 (1) & 4 (2) & \\
Full & 4 (1) & 3 (1) & 2 (1) & 2 (1) & 1 (1) \\
\bottomrule
\end{tabular}
\end{center}
\medskip
\caption{The number of knot types and (in parentheses) the number of possible distinct ropelength values which can be obtained by taking a connect sum of two knots of given symmetry types, assuming that the two  summand knots are not related by a symmetry. For example, the connect sum of a $3_1$ knot (reversible) with a $4_1$ knot (full symmetry) yields two knot types: $3_1 \# 4_1$ and $3_1^m \# 4_1$, but only one ropelength value since these knot types are related by a mirror symmetry. On the other hand, the connect sum of a $3_1$ knot with a $(+)$ amphichiral $8_{17}$ knot yields four possible knot types: $3_1 \# 8_{17}$, $3_1 \# 8_{17}^r$, $3_1^m \# 8_{17}$, $3_1^m \# 8_{17}^r$ but again only one ropelength value since the last three types are related to the first by a reverse, mirror, or mirror-reverse symmetry respectively. On the other hand, the sum of a $3_1$ knot with a (reversible) $5_1$ knot yields four knot types: $3_1 \# 5_1$, $3_1^m \# 5_1$, $3_1 \# 5_1^m$, and $3_1^m \# 5_1^m$ with two potentially different ropelength values, one for $3_1 \# 5_1$ and $3_1^m \# 5_1^m$ (whose tight configurations are related by a mirror symmetry) and one for $3_1^m \# 5_1$ and $3_1 \# 5_1^m$ (where again the tight configurations are related by a mirror symmetry). The tight configurations of the knots $3_1 \# 5_1$ and $3_1^m \# 5_1$ do not seem to be related by a rigid motion and have ropelength values $71.544$ and $71.579$, respectively.}
\label{tab:ktcounts}
\end{table}

We compiled symmetry data for prime knots of 9 and fewer crossings
from Henry/Weeks \cite{MR1164115} and Kodama/Sakuma
\cite{MR1177431}. For a given prime knot ``base'' type, we denoted the
mirror, reversal, and reverse-mirror of the knot (when they are not
isotopic to the base) by the tags $m$, $r$ or $rm$. We then ordered
the list by crossing number, index in the Rolfsen table of
knots~\cite{MR0515288}, and symmetry type, so that $K < K^m < K^r <
K^{rm}$. For composite knots involving two factors, we used the
calculation summarized in Table~\ref{tab:ktcounts} to enumerate the
different knot types possible for the connect sum in terms of the
symmetries of the summands. There were a few cases where we had more
than two summands; these were checked by hand. In our tables, each
composite knot type appears once, labeled with the summands in sorted
order. For example, the label $3_1 \# 3_1^m \# 5_1$ appears in our
list, but the labels $3_1^m \# 3_1 \# 5_1$ and $5_1 \# 3_1 \# 3_1^m$
do not. Mastin~\cite{mastincomposites} provides a general algorithm
for enumerating composites with any number of prime factors and
determining their symmetry types based on the JSJ-decomposition of
composite knots. Tabulating composite \emph{links} is a considerably
more difficult problem, treated in~\cite{cpm_tabulation}.

\subsection{Algorithms}
We minimized polygonal ropelength using the \texttt{ridgerunner} code
described in~\cite{acprtightening}. Recall that we define the
\textit{residual} of an approximately tight polygon to be the fraction
of the gradient of length remaining after projection and that the knot
is critical when this residual vanishes.  For 520 of the 544 tightened
composite knots, the residual values were below $0.01$ (i.e.~over
$99\%$ of the gradient is resolved against the constraints).

Since the minimum length of our composite knots varies considerably
over our table, we did not choose the same number of edges for each
knot. Instead, we choose to keep constant the ``resolution'' of each
polygon, defined as the quotient of the number of edges and the
ropelength of the curve. All of our final configurations have
resolution at least $8$ (from several hundred to around a thousand
vertices).

For each configuration, we give a ropelength upper bound which is
given by carefully numerically approximating the length of a piecewise
$C^2$ curve constructed by splicing short circle arcs into our
polygons. This bound is a rigorous upper bound on the minimum
ropelength of the given knot type, the details of which
appear in~\cite{acprtightening}.

\subsection{Initial configurations} 
The \texttt{ridgerunner} algorithm requires an initial configuration
of each knot type. Since the software proceeds by constrained gradient
descent, it is designed to stop at local minima of the ropelength
function. As a practical matter, it is impossible to know at present
whether any given ropelength critical knot is a global minimum over
its knot type (rigorous, sharp lower bounds are not known for any knot
type, and the configuration space of polygonal curves is far too large
to attempt any kind of exhaustive search). However, we tried to reduce
the probability of false local minima in our dataset in two ways.

First, splicing two given polygons $P_1$ and $P_2$ together requires
a choice of arcs on $P_1$ and $P_2$ to cut out and splice. It is clear
that the shape of the resulting tight composite probably depends on
these choices. To avoid this problem, we took ``all'' connect sums of
$P_1$ and $P_2$ using the following algorithm.

\begin{enumerate}
\item Identify all arcs of $P_1$ and $P_2$ lying on their convex
  hulls, i.e.~lying on the ``outside'' of the configuration.
\item For each pair of arcs, choose an edge on each arc.
\item Translate and rotate the polygons to align the selected edges,
  delete the edges, and splice the polygons together. Repeat this
  process for all pairs of arcs to form an ensemble of composite
  polygons.
\end{enumerate}

For the prime summands, we used the approximately tight configurations
from~\cite{acprtightening}. After splicing the polygons together, we
then smoothed each of these connect sums and scaled them up, allowing
them to retighten from a position with larger thickness. The winning
configurations were selected for further runs at higher resolution.

Second, we explored the configuration space of each of our knots using
a new version of the \texttt{ridgerunner} core called ``mangle
mode''. In this form, rather than attempting to minimize the length of
a polygon subject to the constraints describing the tube, the software
applies a randomly chosen toroidal force field to the knot and
resolves the resulting motion against the tube constraints. This has
the effect of turning the knots ``inside out'', while preserving the tube around the
knot. In practice, using this method to generate twenty alternate
start configurations and minimizing ropelength from each position was
an effective method of discovering alternate local minima for
ropelength.

\subsection{Hardware}
We minimized our knots on the ACCRE cluster at Vanderbilt University and
on a 72 core cluster at the University of St.\ Thomas, running
computations for most of the year 2010 as we experimented with
different start positions and run parameters for the knots. In total,
we ran more than 20,000 configurations, distributed among our
\CompositeKnotCount\ composite knot types.

\section{Results}

We report the ropelengths of our tight shapes in
Tables~\ref{ByKnotTable 0}-\ref{ByKnotTable 5} on
pages~\pageref{ByKnotTable 0}--\pageref{ByKnotTable 5}, and their
residuals in Tables~\ref{ResidualTable 0}--\ref{ResidualTable 5} on
pages~\pageref{ResidualTable 0}--\pageref{ResidualTable 5}.  As
in~\cite{acprtightening}, our data for these composite knots includes
both self-contact sets and measures of the compression force on the
contacts. We hope that these conformations will inspire other groups
to refine them and further investigate them. All of our data,
including the vertices of our approximately tight conformations, are
freely available on Cantarella and Rawdon's web pages.

\subsection{Writhe of composite knots} 
Katritch et al.~\cite{kops} conjectured that the \emph{average writhing
  number} of a composite knot should obey the relation $\Wr(K_1 \#
K_2) = \Wr(K_1) + \Wr(K_2)$. This was a surprising conjecture since
the writhing number generally depends on the entire shape of a
curve. Laing and Sumners~\cite{laingsumners} showed that given two
knots $K_1$ and $K_2$ which intersect in an arc, the conformation of
$K_1 \# K_2$ given by deleting the common arc has writhe equal to
$\Wr(K_1) + \Wr(K_2)$, but there is no guarantee that the \emph{tight}
configuration of $K_1 \# K_2$ should be constructed in such a
way. Nonetheless in the vast majority of the conformations we
computed, it seems that the two prime summands appear almost
unchanged in the tight composite and the writhes do obey this sum
property to a high degree of numerical accuracy.

\begin{figure}
\hfill
\includegraphics[height=1in]{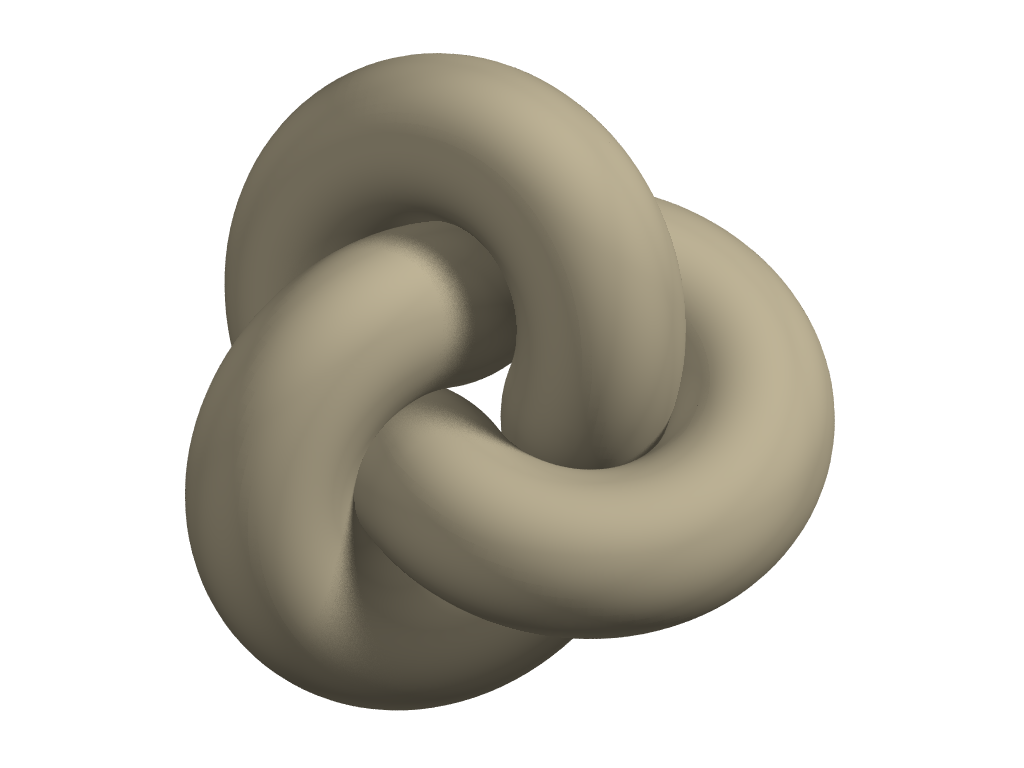} 
\hfill
\includegraphics[height=1.5in]{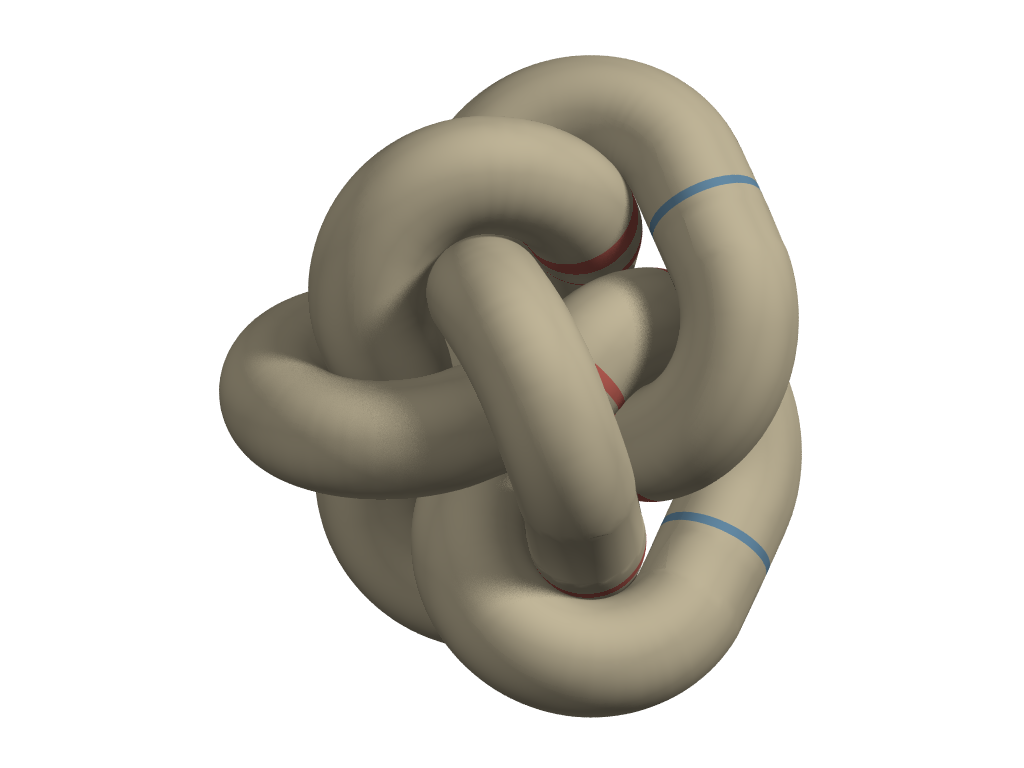} 
\hfill
\includegraphics[height=1.75in]{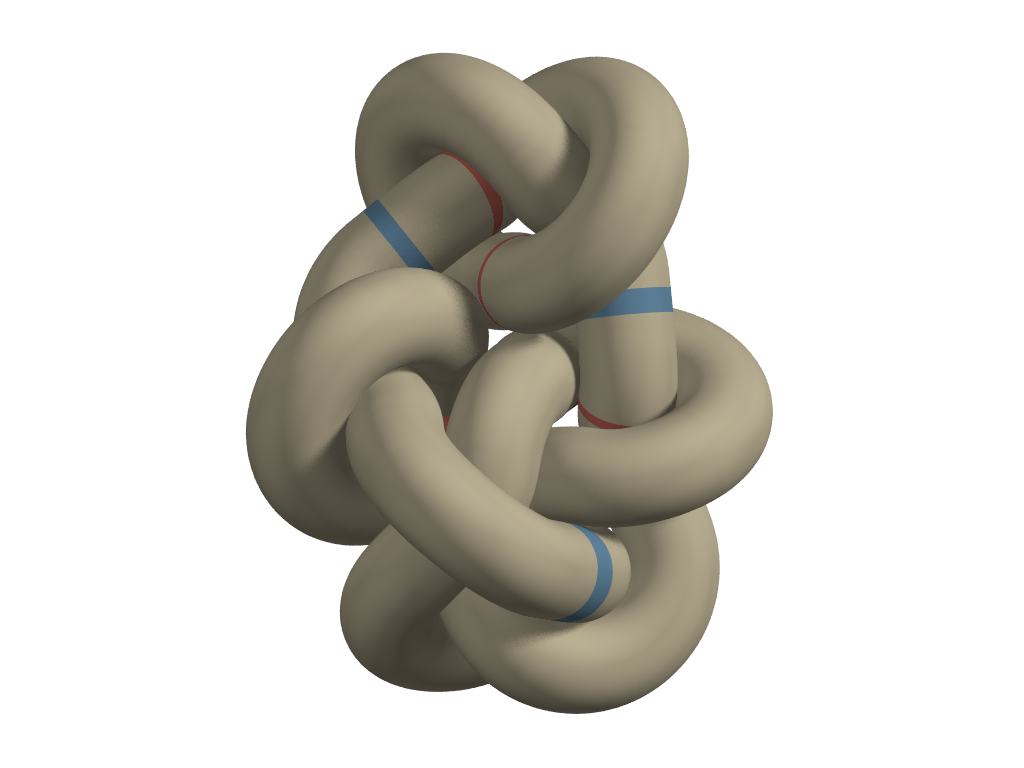}
\hfill
\hphantom{.}
\caption{Tight configurations of $3_1$, $9_{48}^m$
  and their connect sum $3_1 \# 9_{48}^m$. The sum of the writhes of
  the tight $3_1$ and $9_{48}^m$ configurations is $0.523$ while the
  writhe of the tight composite $3_1 \# 9_{48}^m$ is $0.625$. The
  difference is easily explained by looking at the pictures; in the
  connect sum, the lower left arc of the trefoil is spliced to the
  upper right arc of the $9_{48}^m$ shown above. The extra loop of the
  trefoil pushes the remainder of the knot out of its tight shape,
  changing the writhe of the overall composite. The other examples in
  our database look similar; in each case it seems that taking the
  connect sum and then tightening distorts one or both summands. Each
  of these configurations shows kinks (highlighted in red) and
  straight segments (highlighted in blue).}
\label{fig:badwrithe}
\end{figure}

However, in a number of cases, we could not verify the conjecture even
after trying a large number of start configurations. A collection of
these are summarized in Table~\ref{tab:writhe_badguys}. Piotr
Pieranski and Sylwek Przybyl have graciously spot-checked the writhe
and tightening computations in this table using their knot-tightening
and writhe computation codes~\cite{ppcomm}. These checks revealed no
significant difference in writhes, although they were able to tighten
the knots somewhat more using SONO and Przybyl's FEM-based knot
tightener. Although it is possible that we have failed to find the
true minimizer in the cases presented, we think this data
presents a significant challenge to the conjecture and requires
explanation. Figure~\ref{fig:badwrithe} shows a typical example of
this phenomenon.

\begin{table}
\begin{center}
Composite knots which appear to violate the writhe conjecture\\
\begin{tabular}[t]{lrrrc}
\toprule
Knot & $\Wr(K_1 \# K_2)$ & $\Wr(K_1) +\Wr(K_2)$ & Difference (\%) & \# Start Positions\\
\midrule
$ 3_{1} \# 9_{48}^m $ & 0.523 & 0.625 & 19.412 \% & 13 (16,25) \\ 
$ 4_{1} \# 7_{7} $ & -0.534 & -0.632 & 18.359 \% & 55 (25,331) \\ 
$ 3_{1} \# 8_{14}^m $ & 0.556 & 0.634 & 14.127 \% & 50 (16,26) \\ 
$ 3_{1} \# 8_{11} $ & 0.609 & 0.526 & 13.688 \% & 38 (16,21) \\ 
$ 3_{1} \# 6_{2}^m $ & -0.564 & -0.632 & 12.000 \% & 17 (16,27) \\ 
$ 3_{1} \# 8_{16}^m $ & -0.576 & -0.632 & 9.610 \% & 73 (16,23) \\ 
$ 3_{1} \# 9_{45}^m $ & 1.711 & 1.865 & 8.979 \% & 16 (16,27) \\ 
$ 3_{1} \# 9_{22}^m $ & -1.123 & -1.223 & 8.936 \% & 42 (16,24) \\ 
$ 5_{1} \# 7_{2} $ & -0.557 & -0.606 & 8.810 \% & 109 (14,24) \\ 
$ 3_{1} \# 9_{42}^m $ & -2.197 & -2.022 & 7.950 \% & 12 (16,24) \\ 
$ 3_{1} \# 9_{21}^m $ & 1.197 & 1.110 & 7.236 \% & 92 (16,23) \\ 
$ 3_{1} \# 8_{6}^m $ & 0.615 & 0.572 & 6.981 \% & 17 (16,22) \\ 
$ 3_{1} \# 9_{12}^m $ & 1.171 & 1.092 & 6.782 \% & 83 (16,29) \\ 
$ 4_{1} \# 8_{13} $ & -1.116 & -1.189 & 6.542 \% & 105 (25,21) \\ 
$ 5_{2} \# 7_{4}^m $ & 1.170 & 1.234 & 5.420 \% & 21 (27,28) \\ 
$ 3_{1} \# 5_{2}^m $ & 1.183 & 1.134 & 4.200 \% & 16 (16,27) \\ 
$ 5_{2} \# 7_{5}^m $ & 2.771 & 2.886 & 4.130 \% & 98 (27,27) \\ 
\bottomrule
\end{tabular}

\end{center}
\vspace{0.25in}
\caption{This table shows a selection of cases where the writhe conjecture is not supported by our data. The table shows the knot type, writhe of composite, sum of writhes of summands, absolute percentage difference between these numbers, and the number of start configurations tried for the composite (and in parentheses, the two summands). These are not the only cases where we are unable to verify the conjecture, but we have chosen not to show cases with more than two summands (because these more complicated knots simply may not be fully tightened), where the percentage difference was less than $4\%$, and where the writhe of the composite was less than $0.1$ (because in these cases a small absolute difference between writhes close to zero leads to huge percentage differences).}
\label{tab:writhe_badguys}
\end{table}

\subsection{Ropelength of composite knots}
\begin{figure}
\begin{center}
\includegraphics[width=5in]{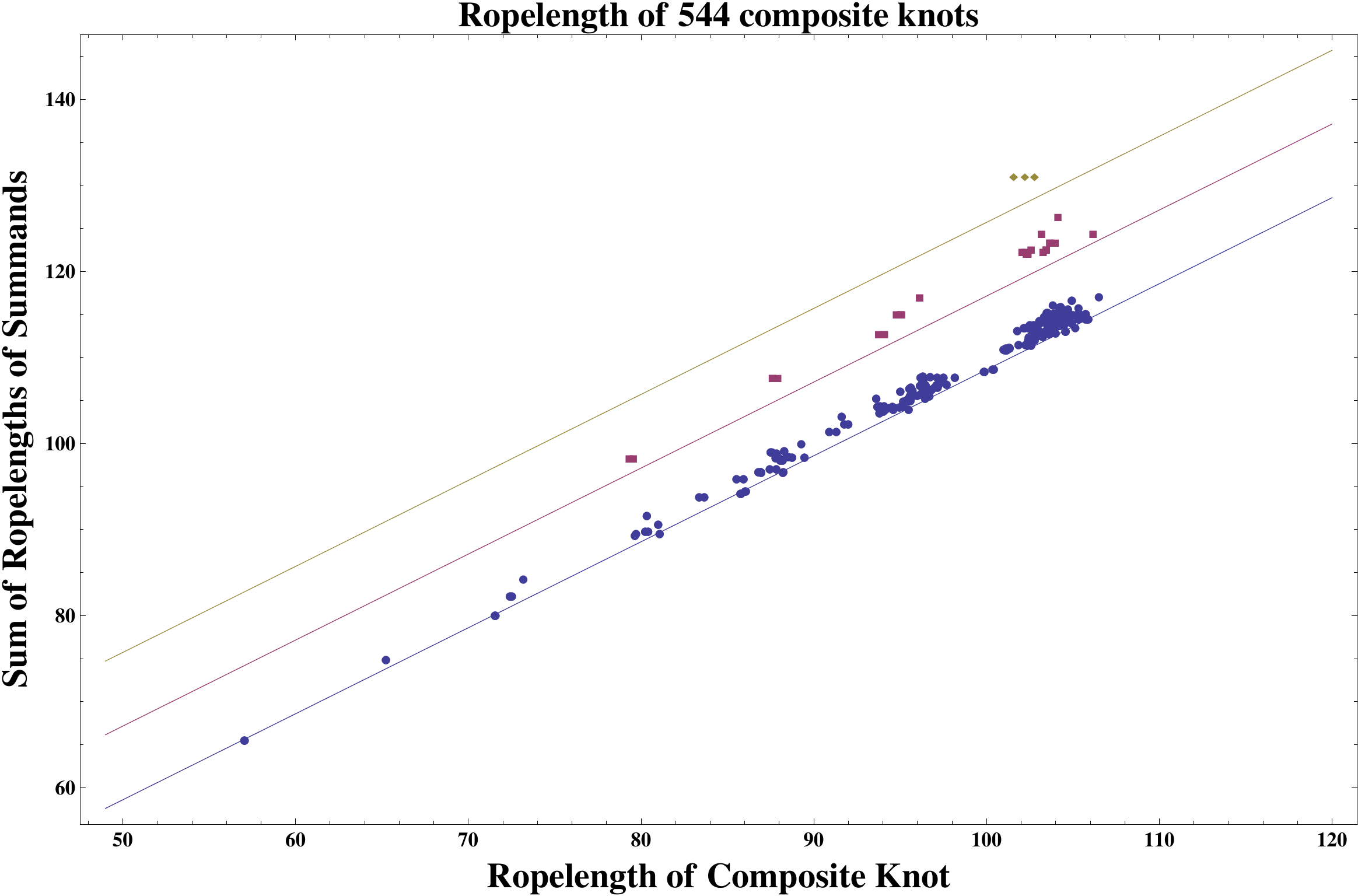}
\end{center}
\caption{This scatterplot shows pairs $(x,y)$ where $x$ is the minimum
  ropelength we have found for a composite knot and $y$ is the sum of
  the minimum ropelengths we have found for its prime
  summands. On the plot, knots with two prime summands are plotted with blue circles, knots with three prime summands are plotted with red squares, and knots with four prime summands are plotted with yellow diamonds. According to the \emph{connect sum conjecture}
  of~\cite{kops}, $y-x$ is at least $(N-1)(4\pi - 4)$, where $N$ is
  the number of prime summands of the knot. This conjecture is shown
  by lines $y-x = 4\pi - 4$, $y-x = 8\pi -8$, and $y-x = 12\pi-12$. As
  one can see from the plot, all our data are very close to obeying
  the bounds predicted by the conjecture.}
\label{fig:compositeropelength}
\end{figure}

Figure~\ref{fig:compositeropelength} shows a scatterplot of the
ropelength of each of our composites ($x$-axis) plotted against the sum
of the ropelengths of their prime summands. Katritch et
al.~\cite{kops} conjectured that the ropelength of the composite
should be at least $4\pi - 4$ less than the sum of the ropelength of
the summands. This has become informally known as the \emph{connect
  sum conjecture} for ropelength. The intuition behind the conjecture
is easy to understand. When two pairs of linked rings (each with
ropelength $8\pi$) are connect-summed to form a three link chain, the
two rings which have been spliced together shrink to form a stadium
curve with ropelength $4\pi + 4$. The difference between the original
ropelength of the rings ($8\pi$) and the ropelength of the stadium
curve ($4\pi + 4$) is the amount of rope saved in the splicing
procedure: $4\pi - 4$. For more complicated knots, it is somewhat surprising that the same amount of rope should be ``exposed'' to a connect sum. If this conjecture holds for very complicated knots, the principle at work would seem to be very different.

However, in this range of composite knots, our data suggests that the conjecture is quite plausible. In our dataset there are only 44 knots where our best conformation of the composite
is very slightly longer than the conjecture predicts. In the most significant example of this phenomenon, the tightest $5_1^m \# 7_1$ has ropelength $100.427$ while the conjecture
predicts that there should be a conformation with ropelength $0.38\%$ lower ($\leq 100.042$). 

\subsection{Other Observations about Tight Configurations}

Table~\ref{tab:bestworst} shows the longest and shortest knots by
crossing number for both prime and composite knots. We can see that
the range of ropelengths for composite knots is smaller than the range
for all knots of the same crossing number. One might expect composite
knots to generally be longer than prime knots of the same crossing
number since composites are separated in two pieces and have less
opportunity to nest together and save rope. But it is mildly
surprising that the longest knots of each crossing number are not
composite. It will be interesting to see if this effect persists
through higher crossing numbers.

\begin{table}
\begin{center}
\begin{tabular}{p{3in}p{3in}}
\begin{center}Composite Knots\end{center} & 
\begin{center}Prime Knots\end{center} \\ \vspace{-0.35in}
\begin{tabular}[t]{lcccc}
\toprule$\Cr$ & $\Rop$ & Links \\
\midrule
$6$ & $ [57.042,57.073] $ & $ 3_{1} \# 3_{1} $, $ 3_{1} \# 3_{1}^m $ \\ 
$7$ & $ [65.240,65.240] $ & $ 3_{1} \# 4_{1} $, $ 3_{1} \# 4_{1} $ \\ 
$8$ & $ [71.544,73.193] $ & $ 3_{1} \# 5_{1} $, $ 4_{1} \# 4_{1} $ \\ 
$9$ & $ [79.329,81.088] $ & $ 3_{1} \# 3_{1} \# 3_{1} $, $ 3_{1} \# 6_{1} $ \\ 
$10$ & $ [85.758,89.472] $ & $ 3_{1} \# 7_{1} $, $ 3_{1} \# 7_{7} $ \\ 
$11$ & $ [83.372,98.171] $ & $ 3_{1} \# 8_{19}^m $, $ 3_{1} \# 8_{18} $ \\ 
$12$ & $ [90.905,106.508] $ & $ 3_{1} \# 9_{46} $, $ 4_{1} \# 8_{18} $ \\ 
\bottomrule
\end{tabular}
 &  \vspace{-0.35in}
\begin{tabular}[t]{lcccc}
\toprule$\Cr$ & $\Rop$ & Knots \\
\midrule
$6$ & $ [56.7058,57.8392] $ & $ 6_{1} $, $ 6_{3} $ \\ 
$7$ & $ [61.4067,65.6924] $ & $ 7_{1} $, $ 7_{6} $ \\ 
$8$ & $ [60.9858,74.9063] $ & $ 8_{19} $, $ 8_{18} $ \\ 
$9$ & $ [68.6169,82.2803] $ & $ 9_{46} $, $ 9_{33} $ \\ 
$10$ & $ [71.0739,92.3565] $ & $ 10_{124} $, $ 10_{123} $ \\ 
$11$ & (no data) & \\
$12$ & (no data) & \\
\bottomrule
\end{tabular}
\end{tabular}
\end{center}
\vspace{0.25in}
\caption{Longest and shortest knots of a given crossing number for prime and composite knots. It is interesting to see that the longest and shortest knots of each crossing number are prime.}
\label{tab:bestworst}
\end{table}

\begin{figure}
\begin{center}
\begin{tabular}{cc}
\includegraphics[height=1.75in]{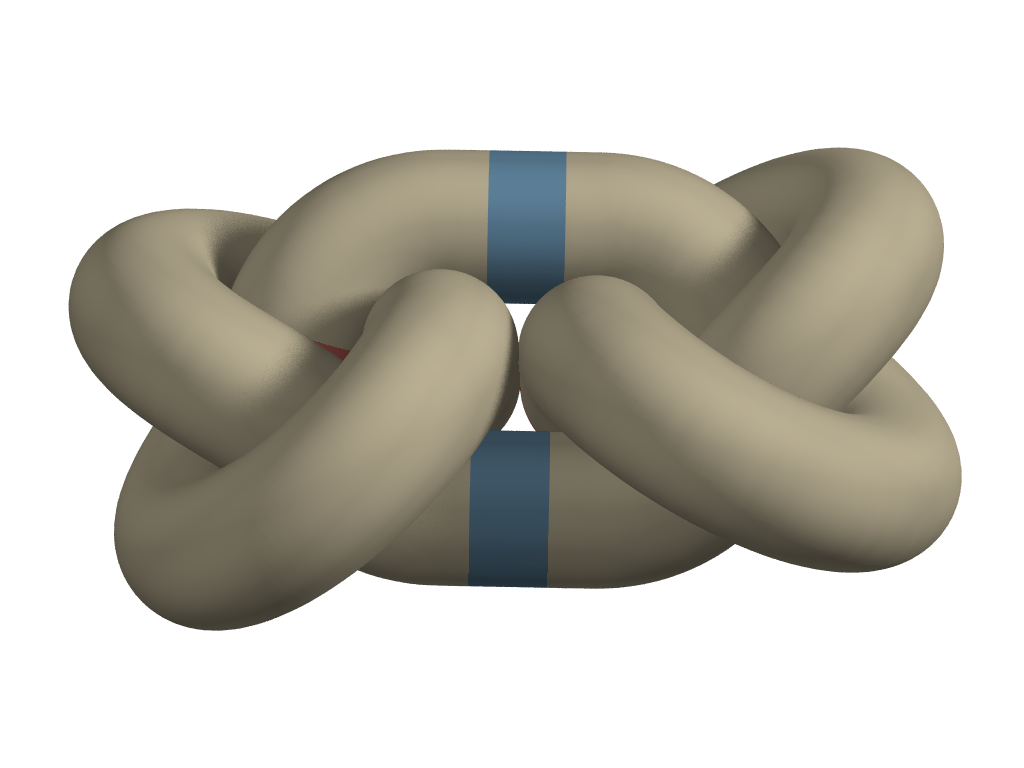} &
\includegraphics[height=1.7in]{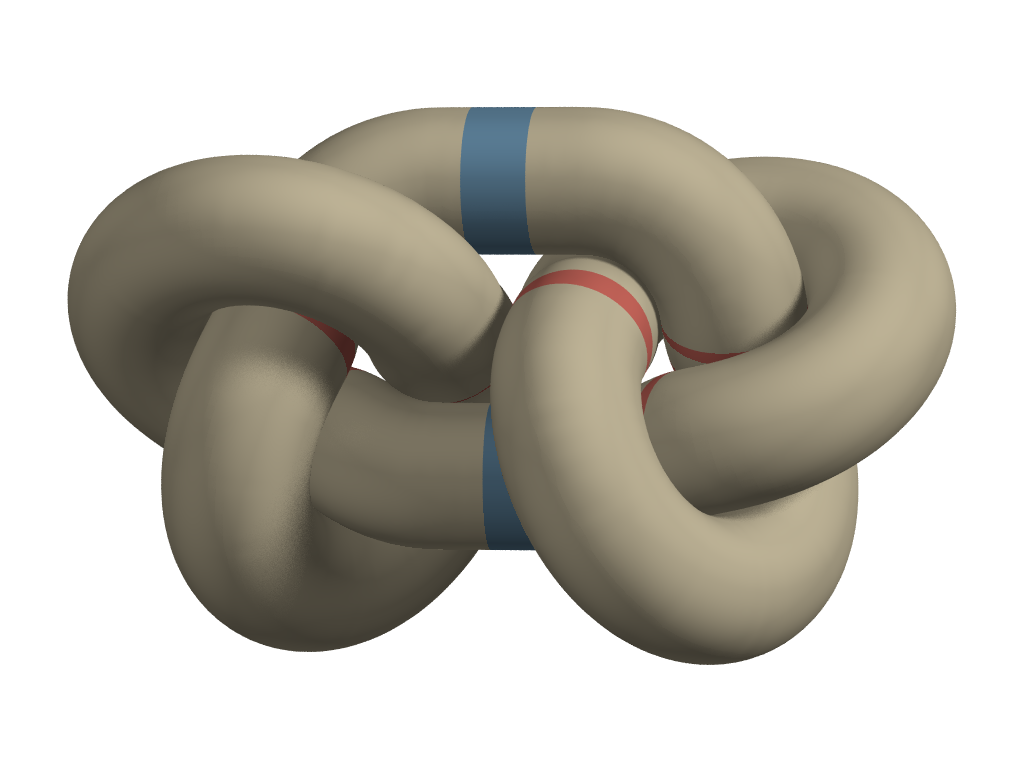}
\end{tabular}
\caption{Two examples of tight composite knots. The left knot is an
  approximately tight configuration of the square knot $3_1 \# 3_1^m$
  which shows the expected straight segments joining the two summands
  highlighted in blue. It is interesting that our tight configuration
  does not have a perfect geometric mirror symmetry (even though a
  critical configuration with this symmetry surely exists;
  see~\cite{bcfm} for a discussion of symmetric criticality). The
  right knot is the approximately right granny knot $3_1 \# 3_1$.  The
  ropelengths are $57.09$ and $57.05$, respectively.}
\label{fig:examples}
\end{center}
\end{figure}

The embedded tube constraint is controlled by both tube-to-tube
contacts (``struts'') and an upper bound on the curvature of the core
polygon (``kinks'')
\cite{acprtightening,MR2001a:57017,MR99k:57025}. It is a theorem
(under some mild regularity assumptions) that no closed
ropelength-critical curve can fail to have a
strut~\cite{2011arXiv1102.3234C}, but kinks seem to be
optional. However, in our data, tight configurations with kinks seem
to be extremely common, occurring in \CompositeKinkCount\ of our
\CompositeKnotCount\ composite configurations (although we expect that
the two configurations without kinks have not fully converged).
Another theorem is that sections of a ropelength-critical curve
without struts or kinks are straight
segments~(see~\cite{2011arXiv1102.3234C,gm} for different versions of
this theorem). Gonzalez conjectured this phenomenon should occur in
all tight composite knots with a mirror symmetry (such as the square
knot in Figure~\ref{fig:examples}). We find this phenomenon in
\CompositeStraightCount\ of the \CompositeKnotCount\ composite knots
with crossing number at most 12.

\section{Future Directions}
Our publicly available dataset of tight knots and links, now including
tight prime knots to 10 crossings, tight prime links to 9 crossings, and (with
this paper) tight composite knots to 12 crossings should provide a
substantial starting point for physicists, biologists, and
mathematicians interested in the geometry of knotted
configurations. We intend to expand the dataset to assemble tight
configurations of all prime and composite knots and links to 12
crossings. This is a substantial undertaking even for prime knots and
links, but the hardest part is certainly composite links. The problem
is that composite links remain untabulated (it is still a challenging
open problem to come up with an algorithm for tabulating composite
links; see~\cite{cpm_tabulation}).

It is certainly interesting that so many of our configurations exhibit
``kinked'' sections of maximum curvature. An excellent confirmation of
this phenomenon would be to rerun our configurations with the
curvature constraint removed to see whether we can achieve shorter
lengths. In addition, the FEM techniques of Pieranski and Przybyl show
great promise in computing at very high resolutions (up to tens of
thousands of vertices). We intend to collaborate with this group for
our next set of computations, using \texttt{ridgerunner} as a
medium-resolution search tool to explore the configuration space of
curves of a given knot type and then switching to FEM for the final
tightening.

\noindent Acknowledgements. We are grateful to all the members of the UGA
Geometry VIGRE group who contributed their time and effort towards
this project. In particular, we would like to acknowledge the
contributions of graduate students Ted Ashton and Steve Lane, who
compiled the original table of composite knots, and Michael Piatek,
who cowrote \texttt{ridgerunner}. We would also like to thank Piotr
Pieranski and Sylwek Przybyl for checking our calculations. Our work
was supported by the UGA VIGRE grants DMS-07-38586 and DMS-00-89927,
UGA REU site grant DMS-06-49242, and NSF grants DMS-06-21903,
DMS-08-10415, and DMS-11-15722 (to Rawdon).

\newpage

\section*{References}
\bibliography{rr05siggraph,drl,cantarella,tsnnls,wonderbib} 

\begin{thebibliography}{10}

\bibitem{MR99b:92032}
Vsevolod Katritch, Jan Bednar, Didier Michoud, Robert~G. Scharein, Jacques
  Dubochet, and Andrzej Stasiak.
\newblock Geometry and physics of knots.
\newblock {\em Nature}, 384(6605):142--145, 1996.

\bibitem{mbms}
Cristian Micheletti, Jaynath Banavar, Amos Maritan, and F.~Seno.
\newblock Protein structures and optimal folding from a geometric variational
  principle.
\newblock {\em Physical Review Letters}, 82:3372--3375, 1999.

\bibitem{MR1981020}
Akos Dobay, Jacques Dubochet, Kenneth Millett, Pierre-Edouard Sottas, and
  Andrzej Stasiak.
\newblock Scaling behavior of random knots.
\newblock {\em Proc. Natl. Acad. Sci. USA}, 100(10):5611--5615 (electronic),
  2003.

\bibitem{1367-2630-3-1-310}
Piotr Pieranski, Sandor Kasas, Giovanni Dietler, Jacques Dubochet, and Andrzej
  Stasiak.
\newblock Localization of breakage points in knotted strings.
\newblock {\em New Journal of Physics}, 3:10, 2001.

\bibitem{kepglue}
Roman~V. Buniy and Thomas~W. Kephart.
\newblock A model of glueballs.
\newblock {\em Phys. Lett.}, B576:127--134, 2003.

\bibitem{MR2003h:58014}
Jason Cantarella, Robert~B. Kusner, and John~M. Sullivan.
\newblock On the minimum ropelength of knots and links.
\newblock {\em Invent. Math.}, 150(2):257--286, 2002.

\bibitem{rawdonthesis}
Eric Rawdon.
\newblock {\em The Thickness of Polygonal Knots}.
\newblock PhD thesis, The University of Iowa, 1997.

\bibitem{MR1749504}
Eric~J. Rawdon.
\newblock Approximating smooth thickness.
\newblock {\em J. Knot Theory Ramifications}, 9(1):113--145, 2000.

\bibitem{MR2207788}
Elizabeth Denne, Yuanan Diao, and John~M. Sullivan.
\newblock Quadrisecants give new lower bounds for the ropelength of a knot.
\newblock {\em Geom. Topol.}, 10:1--26 (electronic), 2006.

\bibitem{MR2471403}
J.~Baranska, S.~Przybyl, and P.~Pieranski.
\newblock Curvature and torsion of the tight closed trefoil knot.
\newblock {\em Eur. Phys. J. B}, 66(4):547--556, 2008.

\bibitem{MR2034393}
Eric~J. Rawdon.
\newblock Can computers discover ideal knots?
\newblock {\em Experiment. Math.}, 12(3):287--302, 2003.

\bibitem{MR1702020}
Andrzej Stasiak, Jacques Dubochet, Vsevolod Katritch, and Piotr Pieranski.
\newblock Ideal knots and their relation to the physics of real knots.
\newblock In {\em Ideal knots}, volume~19 of {\em Ser. Knots Everything}, pages
  1--19. World Sci. Publishing, River Edge, NJ, 1998.

\bibitem{tightopen}
Piotr Pieranski, Sylwester Przybyl, and Eric Rawdon.
\newblock Tight open knots.
\newblock {\em The European Physical Journal E - Soft Matter}, 6:123--128,
  2001.

\bibitem{baranska}
Justyna Baranska, Piotr Pieranski, Sylwester Przybyl, and Eric~J. Rawdon.
\newblock Length of the tightest trefoil knot.
\newblock {\em Physical Review E (Statistical, Nonlinear, and Soft Matter
  Physics)}, 70(5):051810, 2004.

\bibitem{gordian}
P.~Pieranski, S.~Przybyl, and A.~Stasiak.
\newblock Gordian unknots.
\newblock arXiv: physics/0103080, 2004.

\bibitem{smutnyphd}
Jana Smutny.
\newblock {\em Global Radii of Curvature, and the Biarc Approximation of Space
  Curves: In Pursuit of Ideal Knot Shapes}.
\newblock PhD thesis, Ecole Polytechnique F\'ed\'erale de Lausanne, 2004.

\bibitem{maddocksbook}
Mathias Carlen, Ben Laurie, John~H. Maddocks, and Jana Smutny.
\newblock Biarcs, global radius of curvature, and the computation of ideal knot
  shapes.
\newblock In {\em Physical and Numerical Models in Knot Theory}, volume~36 of
  {\em Ser. Knots Everything}, pages 75--108. World Sci. Publishing, Singapore,
  2005.

\bibitem{acprtightening}
Ted Ashton, Jason Cantarella, Michael Piatek, and Eric~J. Rawdon.
\newblock Knot tightening by constrained gradient descent.
\newblock {\em Experimental Mathematics}, 20(1):57--90, 2011.

\bibitem{tait}
Peter~Guthrie Tait.
\newblock Sevenfold knottiness.
\newblock {\em Proc. Royal Soc. Edinburgh}, 9(98):363--366, 1876-7.

\bibitem{kirkman}
Kirkman.
\newblock The enumeration, description, and construction of knots of fewer than
  10 crossings.
\newblock {\em Trans. Roy. Soc. Edinburgh}, 32:281--309, 1883-4.

\bibitem{MR2105106}
Slavik~V. Jablan.
\newblock New knot tables.
\newblock {\em Filomat}, (15):141--152, 2001.

\bibitem{MR1646740}
Jim Hoste, Morwen Thistlethwaite, and Jeff Weeks.
\newblock The first 1,701,936 knots.
\newblock {\em Math. Intelligencer}, 20(4):33--48, 1998.

\bibitem{kops}
Vsevolod Katritch, Wilma Olson, Piotr Pieranski, Jacques Dubochet, and Andrzej
  Stasiak.
\newblock Properties of ideal composite knots.
\newblock {\em Nature}, 388(6638):148--151, 1997.

\bibitem{MR0031733}
Horst Schubert.
\newblock Die eindeutige {Z}erlegbarkeit eines {K}notens in {P}rimknoten.
\newblock {\em S.-B. Heidelberger Akad. Wiss. Math.-Nat. Kl.}, 1949(3):57--104,
  1949.

\bibitem{MR0258014}
J.~H. Conway.
\newblock An enumeration of knots and links, and some of their algebraic
  properties.
\newblock In {\em Computational {P}roblems in {A}bstract {A}lgebra ({P}roc.
  {C}onf., {O}xford, 1967)}, pages 329--358. Pergamon, Oxford, 1970.

\bibitem{MR1164115}
Shawn~R. Henry and Jeffrey~R. Weeks.
\newblock Symmetry groups of hyperbolic knots and links.
\newblock {\em J. Knot Theory Ramifications}, 1(2):185--201, 1992.

\bibitem{MR1177431}
Kouzi Kodama and Makoto Sakuma.
\newblock Symmetry groups of prime knots up to {$10$} crossings.
\newblock In {\em Knots 90 ({O}saka, 1990)}, pages 323--340. de Gruyter,
  Berlin, 1992.

\bibitem{MR0515288}
Dale Rolfsen.
\newblock {\em Knots and links}.
\newblock Publish or Perish Inc., Berkeley, Calif., 1976.
\newblock Mathematics Lecture Series, No. 7.

\bibitem{mastincomposites}
Matt Mastin.
\newblock Composite knots and their symmetries.
\newblock In preparation.

\bibitem{cpm_tabulation}
Jason Cantarella, Jason Parsley, and Matt Mastin.
\newblock Symmetries and tabulation of composite links.
\newblock In preparation.

\bibitem{laingsumners}
Christian Laing, Renzo Ricca, and De~Witt Sumners.
\newblock Writhe, helicity and vortex reconnection.
\newblock in preparation.

\bibitem{ppcomm}
Piotr Pieranski and Sylwek Przybyl.
\newblock personal communication.

\bibitem{MR2001a:57017}
Eric~J. Rawdon.
\newblock Approximating smooth thickness.
\newblock {\em J. Knot Theory Ramifications}, 9(1):113--145, 2000.

\bibitem{MR99k:57025}
R.~A. Litherland, J.~Simon, O.~Durumeric, and E.~Rawdon.
\newblock Thickness of knots.
\newblock {\em Topology Appl.}, 91(3):233--244, 1999.

\bibitem{2011arXiv1102.3234C}
J.~{Cantarella}, J.~H.~G. {Fu}, R.~{Kusner}, and J.~M. {Sullivan}.
\newblock {Ropelength Criticality}.
\newblock {\em ArXiv e-prints}, February 2011.

\bibitem{gm}
Oscar Gonzalez and John~H. Maddocks.
\newblock Global curvature, thickness, and the ideal shapes of knots.
\newblock {\em Proc. Nat. Acad. Sci. (USA)}, 96:4769--4773, 1999.

\bibitem{bcfm}
Jennifer {Belton}, Jason {Cantarella}, Joseph~H.G. {Fu}, , and Matt {Mastin}.
\newblock {Symmetric Criticality for Tight Knots}.
\newblock In preparation.

\end{thebibliography}

\appendix
\section*{Ropelength Tables}
The pages that follow contain two sets of tables of ropelength
data. The first set, Tables~\ref{ByKnotTable 0}-\ref{ByKnotTable 5} on
pages~\pageref{ByKnotTable 0}--\pageref{ByKnotTable 5}, show the
polygonal ropelength ($\PRop$) and ropelength upper bounds ($\Rop$)
that we have obtained for each of the composite knot types that we
have considered. The composite knots are organized by dictionary order
on their summands, with the primary order coming from position in
Rolfsen's table \cite{MR0515288}, with the knot $X^y_z$ being the $z$-th example of a
prime $X$-crossing link of $y$ components in the table. Knots with the
same Rolfsen position are ordered by symmetry type, with the
convention $K < K^m < K^r < K^{rm}$. To aid the reader
in making sense of the table, we insert lines where the crossing
number changes and spaces where the base types of the summands change.

The second set of tables, Tables~\ref{ResidualTable
  0}--\ref{ResidualTable 5} on pages~\pageref{ResidualTable
  0}--\pageref{ResidualTable 5} give the residual of each of our
computed configurations. The low residuals show that they are close to
critical. We include this data as a measure of the relative quality of
each of our minimized configurations.

\begin{table}[h]
 \caption{\label{ByKnotTable 0} Ropelengths of Tight Knots  by Knot Type, Part 1 of 6}

\begin{center}
\hfill
\begin{tabular}[t]{l@{\hspace{0.1in}}r@{\hspace{0.1in}}r} \toprule
Knot & \multicolumn{1}{l}{$\PRop$} & \multicolumn{1}{l}{$\Rop$}\\ \midrule
$3_{1} \# 3_{1}$ & $57.05$ & $57.04$ \\
 
$3_{1} \# 3_{1}^m$ & $57.09$ & $57.07$ \\
 
$3_{1}^m \# 3_{1}^m$ & $57.05$ & $57.04$ \\
\addlinespace[0.589183333333334em]\cmidrule(r{0.5em}l{0.5em}){1-3}\addlinespace[0.589183333333334em]

$3_{1} \# 4_{1}$ & $65.26$ & $65.24$ \\
 
$3_{1}^m \# 4_{1}$ & $65.26$ & $65.24$ \\
\addlinespace[0.589183333333334em]\cmidrule(r{0.5em}l{0.5em}){1-3}\addlinespace[0.589183333333334em]

$3_{1} \# 5_{1}$ & $71.56$ & $71.54$ \\
 
$3_{1} \# 5_{1}^m$ & $71.59$ & $71.58$ \\
 
$3_{1}^m \# 5_{1}$ & $71.59$ & $71.58$ \\
 
$3_{1}^m \# 5_{1}^m$ & $71.56$ & $71.54$ \\
\addlinespace[0.589183333333334em]

$3_{1} \# 5_{2}$ & $72.43$ & $72.41$ \\
 
$3_{1} \# 5_{2}^m$ & $72.55$ & $72.53$ \\
 
$3_{1}^m \# 5_{2}$ & $72.55$ & $72.53$ \\
 
$3_{1}^m \# 5_{2}^m$ & $72.43$ & $72.41$ \\
\addlinespace[0.589183333333334em]

$4_{1} \# 4_{1}$ & $73.20$ & $73.19$ \\
\addlinespace[0.589183333333334em]\cmidrule(r{0.5em}l{0.5em}){1-3}\addlinespace[0.589183333333334em]

$3_{1} \# 6_{1}$ & $81.10$ & $81.09$ \\
 
$3_{1} \# 6_{1}^m$ & $79.73$ & $79.72$ \\
 
$3_{1}^m \# 6_{1}$ & $79.73$ & $79.72$ \\
 
$3_{1}^m \# 6_{1}^m$ & $81.10$ & $81.09$ \\
\addlinespace[0.589183333333334em]

$3_{1} \# 6_{2}$ & $80.27$ & $80.25$ \\
 
$3_{1} \# 6_{2}^m$ & $80.43$ & $80.42$ \\
 
$3_{1}^m \# 6_{2}$ & $80.43$ & $80.42$ \\
 
$3_{1}^m \# 6_{2}^m$ & $80.27$ & $80.25$ \\
\addlinespace[0.589183333333334em]

$3_{1} \# 6_{3}$ & $81.02$ & $81.00$ \\
 
$3_{1}^m \# 6_{3}$ & $81.02$ & $81.00$ \\
\addlinespace[0.589183333333334em]

$4_{1} \# 5_{1}$ & $79.67$ & $79.65$ \\
 
$4_{1} \# 5_{1}^m$ & $79.67$ & $79.65$ \\
\addlinespace[0.589183333333334em]

$4_{1} \# 5_{2}$ & $80.36$ & $80.34$ \\
 
$4_{1} \# 5_{2}^m$ & $80.36$ & $80.34$ \\

\bottomrule
\end{tabular}
\hfill\begin{tabular}[t]{l@{\hspace{0.1in}}r@{\hspace{0.1in}}r} \toprule
Knot & $\PRop$ & $\Rop$\\ \midrule
$3_{1} \# 3_{1} \# 3_{1}$ & $79.35$ & $79.33$ \\
 
$3_{1} \# 3_{1} \# 3_{1}^m$ & $79.57$ & $79.56$ \\
 
$3_{1} \# 3_{1}^m \# 3_{1}^m$ & $79.57$ & $79.56$ \\
 
$3_{1}^m \# 3_{1}^m \# 3_{1}^m$ & $79.35$ & $79.33$ \\
\addlinespace[0.461081041666666em]\cmidrule(r{0.5em}l{0.5em}){1-3}\addlinespace[0.461081041666666em]

$3_{1} \# 7_{1}$ & $85.78$ & $85.76$ \\
 
$3_{1} \# 7_{1}^m$ & $85.82$ & $85.80$ \\
 
$3_{1}^m \# 7_{1}$ & $85.82$ & $85.80$ \\
 
$3_{1}^m \# 7_{1}^m$ & $85.78$ & $85.76$ \\
\addlinespace[0.461081041666666em]

$3_{1} \# 7_{2}$ & $86.98$ & $86.96$ \\
 
$3_{1} \# 7_{2}^m$ & $88.24$ & $88.22$ \\
 
$3_{1}^m \# 7_{2}$ & $88.24$ & $88.22$ \\
 
$3_{1}^m \# 7_{2}^m$ & $86.98$ & $86.96$ \\
\addlinespace[0.461081041666666em]

$3_{1} \# 7_{3}$ & $88.27$ & $88.25$ \\
 
$3_{1} \# 7_{3}^m$ & $86.94$ & $86.92$ \\
 
$3_{1}^m \# 7_{3}$ & $86.94$ & $86.92$ \\
 
$3_{1}^m \# 7_{3}^m$ & $88.27$ & $88.25$ \\
\addlinespace[0.461081041666666em]

$3_{1} \# 7_{4}$ & $87.85$ & $87.84$ \\
 
$3_{1} \# 7_{4}^m$ & $87.48$ & $87.46$ \\
 
$3_{1}^m \# 7_{4}$ & $87.48$ & $87.46$ \\
 
$3_{1}^m \# 7_{4}^m$ & $87.85$ & $87.84$ \\
\addlinespace[0.461081041666666em]

$3_{1} \# 7_{5}$ & $88.22$ & $88.20$ \\
 
$3_{1} \# 7_{5}^m$ & $88.09$ & $88.07$ \\
 
$3_{1}^m \# 7_{5}$ & $88.09$ & $88.07$ \\
 
$3_{1}^m \# 7_{5}^m$ & $88.22$ & $88.20$ \\
\addlinespace[0.461081041666666em]

$3_{1} \# 7_{6}$ & $88.57$ & $88.55$ \\
 
$3_{1} \# 7_{6}^m$ & $88.45$ & $88.43$ \\
 
$3_{1}^m \# 7_{6}$ & $88.45$ & $88.43$ \\
 
$3_{1}^m \# 7_{6}^m$ & $88.57$ & $88.55$ \\
\addlinespace[0.461081041666666em]

$3_{1} \# 7_{7}$ & $89.49$ & $89.47$ \\
 
$3_{1} \# 7_{7}^m$ & $88.77$ & $88.75$ \\
 
$3_{1}^m \# 7_{7}$ & $88.77$ & $88.75$ \\

\bottomrule
\end{tabular}
\hfill\begin{tabular}[t]{l@{\hspace{0.1in}}r@{\hspace{0.1in}}r} \toprule
Knot & $\PRop$ & $\Rop$\\ \midrule
$3_{1}^m \# 7_{7}^m$ & $89.49$ & $89.47$ \\
\addlinespace[0.503841527777778em]

$4_{1} \# 6_{1}$ & $87.88$ & $87.86$ \\
 
$4_{1} \# 6_{1}^m$ & $87.88$ & $87.86$ \\
\addlinespace[0.503841527777778em]

$4_{1} \# 6_{2}$ & $88.31$ & $88.29$ \\
 
$4_{1} \# 6_{2}^m$ & $88.31$ & $88.29$ \\
\addlinespace[0.503841527777778em]

$4_{1} \# 6_{3}$ & $89.31$ & $89.29$ \\
\addlinespace[0.503841527777778em]

$5_{1} \# 5_{1}$ & $86.06$ & $86.03$ \\
 
$5_{1} \# 5_{1}^m$ & $86.11$ & $86.09$ \\
 
$5_{1}^m \# 5_{1}^m$ & $86.06$ & $86.03$ \\
\addlinespace[0.503841527777778em]

$5_{1} \# 5_{2}$ & $86.82$ & $86.80$ \\
 
$5_{1} \# 5_{2}^m$ & $86.97$ & $86.95$ \\
 
$5_{1}^m \# 5_{2}$ & $86.97$ & $86.95$ \\
 
$5_{1}^m \# 5_{2}^m$ & $86.82$ & $86.80$ \\
\addlinespace[0.503841527777778em]

$5_{2} \# 5_{2}$ & $87.52$ & $87.51$ \\
 
$5_{2} \# 5_{2}^m$ & $87.61$ & $87.60$ \\
 
$5_{2}^m \# 5_{2}^m$ & $87.52$ & $87.51$ \\
\addlinespace[0.503841527777778em]\cmidrule(r{1.5em}l{1.5em}){1-3}\addlinespace[0.503841527777778em]

$3_{1} \# 3_{1} \# 4_{1}$ & $87.63$ & $87.61$ \\
 
$3_{1} \# 3_{1}^m \# 4_{1}$ & $87.94$ & $87.92$ \\
 
$3_{1}^m \# 3_{1}^m \# 4_{1}$ & $87.63$ & $87.61$ \\
\addlinespace[0.503841527777778em]\cmidrule(r{0.5em}l{0.5em}){1-3}\addlinespace[0.503841527777778em]

$3_{1} \# 8_{1}$ & $93.92$ & $93.90$ \\
 
$3_{1} \# 8_{1}^m$ & $94.06$ & $94.04$ \\
 
$3_{1}^m \# 8_{1}$ & $94.06$ & $94.04$ \\
 
$3_{1}^m \# 8_{1}^m$ & $93.92$ & $93.90$ \\
\addlinespace[0.503841527777778em]

$3_{1} \# 8_{2}$ & $94.42$ & $94.40$ \\
 
$3_{1} \# 8_{2}^m$ & $94.68$ & $94.66$ \\
 
$3_{1}^m \# 8_{2}$ & $94.68$ & $94.66$ \\
 
$3_{1}^m \# 8_{2}^m$ & $94.42$ & $94.40$ \\
\addlinespace[0.503841527777778em]

$3_{1} \# 8_{3}$ & $95.53$ & $95.50$ \\
 
$3_{1}^m \# 8_{3}$ & $95.53$ & $95.50$ \\

\bottomrule
\end{tabular}
\hfill\hfill\end{center}
\end{table}

\newpage
\begin{table}[h]
 \caption{\label{ByKnotTable 1} Ropelengths of Tight Knots  by Knot Type, Part 2 of 6}

\begin{center}
\hfill\begin{tabular}[t]{l@{\hspace{0.1in}}r@{\hspace{0.1in}}r} \toprule
Knot & $\PRop$ & $\Rop$\\ \midrule
$3_{1} \# 8_{4}$ & $94.63$ & $94.61$ \\
 
$3_{1} \# 8_{4}^m$ & $94.98$ & $94.96$ \\
 
$3_{1}^m \# 8_{4}$ & $94.98$ & $94.96$ \\
 
$3_{1}^m \# 8_{4}^m$ & $94.63$ & $94.61$ \\
\addlinespace[0.499734999999999em]

$3_{1} \# 8_{5}$ & $95.20$ & $95.18$ \\
 
$3_{1} \# 8_{5}^m$ & $95.40$ & $95.38$ \\
 
$3_{1}^m \# 8_{5}$ & $95.40$ & $95.38$ \\
 
$3_{1}^m \# 8_{5}^m$ & $95.20$ & $95.18$ \\
\addlinespace[0.499734999999999em]

$3_{1} \# 8_{6}$ & $95.62$ & $95.59$ \\
 
$3_{1} \# 8_{6}^m$ & $95.64$ & $95.62$ \\
 
$3_{1}^m \# 8_{6}$ & $95.64$ & $95.62$ \\
 
$3_{1}^m \# 8_{6}^m$ & $95.62$ & $95.59$ \\
\addlinespace[0.499734999999999em]

$3_{1} \# 8_{7}$ & $95.48$ & $95.47$ \\
 
$3_{1} \# 8_{7}^m$ & $95.62$ & $95.59$ \\
 
$3_{1}^m \# 8_{7}$ & $95.62$ & $95.59$ \\
 
$3_{1}^m \# 8_{7}^m$ & $95.48$ & $95.47$ \\
\addlinespace[0.499734999999999em]

$3_{1} \# 8_{8}$ & $96.94$ & $96.69$ \\
 
$3_{1} \# 8_{8}^m$ & $95.58$ & $95.56$ \\
 
$3_{1}^m \# 8_{8}$ & $95.58$ & $95.56$ \\
 
$3_{1}^m \# 8_{8}^m$ & $96.94$ & $96.69$ \\
\addlinespace[0.499734999999999em]

$3_{1} \# 8_{9}$ & $96.47$ & $96.45$ \\
 
$3_{1}^m \# 8_{9}$ & $96.47$ & $96.45$ \\
\addlinespace[0.499734999999999em]

$3_{1} \# 8_{10}$ & $96.22$ & $96.20$ \\
 
$3_{1} \# 8_{10}^m$ & $96.65$ & $96.63$ \\
 
$3_{1}^m \# 8_{10}$ & $96.65$ & $96.63$ \\
 
$3_{1}^m \# 8_{10}^m$ & $96.22$ & $96.20$ \\
\addlinespace[0.499734999999999em]

$3_{1} \# 8_{11}$ & $95.70$ & $95.68$ \\
 
$3_{1} \# 8_{11}^m$ & $95.64$ & $95.62$ \\
 
$3_{1}^m \# 8_{11}$ & $95.64$ & $95.62$ \\
 
$3_{1}^m \# 8_{11}^m$ & $95.70$ & $95.68$ \\
\addlinespace[0.499734999999999em]

$3_{1} \# 8_{12}$ & $96.48$ & $96.47$ \\

\bottomrule
\end{tabular}
\hfill\begin{tabular}[t]{l@{\hspace{0.1in}}r@{\hspace{0.1in}}r} \toprule
Knot & $\PRop$ & $\Rop$\\ \midrule
$3_{1}^m \# 8_{12}$ & $96.48$ & $96.47$ \\
\addlinespace[0.499734999999999em]

$3_{1} \# 8_{13}$ & $96.44$ & $96.42$ \\
 
$3_{1} \# 8_{13}^m$ & $96.02$ & $95.99$ \\
 
$3_{1}^m \# 8_{13}$ & $96.02$ & $95.99$ \\
 
$3_{1}^m \# 8_{13}^m$ & $96.44$ & $96.42$ \\
\addlinespace[0.499734999999999em]

$3_{1} \# 8_{14}$ & $97.19$ & $97.17$ \\
 
$3_{1} \# 8_{14}^m$ & $97.03$ & $97.00$ \\
 
$3_{1}^m \# 8_{14}$ & $97.03$ & $97.00$ \\
 
$3_{1}^m \# 8_{14}^m$ & $97.19$ & $97.17$ \\
\addlinespace[0.499734999999999em]

$3_{1} \# 8_{15}$ & $96.54$ & $96.52$ \\
 
$3_{1} \# 8_{15}^m$ & $96.33$ & $96.31$ \\
 
$3_{1}^m \# 8_{15}$ & $96.33$ & $96.31$ \\
 
$3_{1}^m \# 8_{15}^m$ & $96.54$ & $96.52$ \\
\addlinespace[0.499734999999999em]

$3_{1} \# 8_{16}$ & $96.51$ & $96.50$ \\
 
$3_{1} \# 8_{16}^m$ & $96.84$ & $96.82$ \\
 
$3_{1}^m \# 8_{16}$ & $96.84$ & $96.82$ \\
 
$3_{1}^m \# 8_{16}^m$ & $96.51$ & $96.50$ \\
\addlinespace[0.499734999999999em]

$3_{1} \# 8_{17}$ & $97.44$ & $97.41$ \\
 
$3_{1} \# 8_{17}^r$ & $97.44$ & $97.41$ \\
 
$3_{1}^m \# 8_{17}$ & $97.44$ & $97.41$ \\
 
$3_{1}^m \# 8_{17}^r$ & $97.44$ & $97.41$ \\
\addlinespace[0.499734999999999em]

$3_{1} \# 8_{18}$ & $98.19$ & $98.17$ \\
 
$3_{1}^m \# 8_{18}$ & $98.19$ & $98.17$ \\
\addlinespace[0.499734999999999em]

$3_{1} \# 8_{19}$ & $83.68$ & $83.67$ \\
 
$3_{1} \# 8_{19}^m$ & $83.39$ & $83.37$ \\
 
$3_{1}^m \# 8_{19}$ & $83.39$ & $83.37$ \\
 
$3_{1}^m \# 8_{19}^m$ & $83.68$ & $83.67$ \\
\addlinespace[0.499734999999999em]

$3_{1} \# 8_{20}$ & $85.55$ & $85.53$ \\
 
$3_{1} \# 8_{20}^m$ & $85.95$ & $85.94$ \\
 
$3_{1}^m \# 8_{20}$ & $85.95$ & $85.94$ \\
 
$3_{1}^m \# 8_{20}^m$ & $85.55$ & $85.53$ \\

\bottomrule
\end{tabular}
\hfill\begin{tabular}[t]{l@{\hspace{0.1in}}r@{\hspace{0.1in}}r} \toprule
Knot & $\PRop$ & $\Rop$\\ \midrule
$3_{1} \# 8_{21}$ & $87.81$ & $87.79$ \\
 
$3_{1} \# 8_{21}^m$ & $88.01$ & $87.99$ \\
 
$3_{1}^m \# 8_{21}$ & $88.01$ & $87.99$ \\
 
$3_{1}^m \# 8_{21}^m$ & $87.81$ & $87.79$ \\
\addlinespace[0.484656666666666em]

$4_{1} \# 7_{1}$ & $93.83$ & $93.81$ \\
 
$4_{1} \# 7_{1}^m$ & $93.83$ & $93.81$ \\
\addlinespace[0.484656666666666em]

$4_{1} \# 7_{2}$ & $96.31$ & $96.29$ \\
 
$4_{1} \# 7_{2}^m$ & $96.31$ & $96.29$ \\
\addlinespace[0.484656666666666em]

$4_{1} \# 7_{3}$ & $95.04$ & $95.02$ \\
 
$4_{1} \# 7_{3}^m$ & $95.04$ & $95.02$ \\
\addlinespace[0.484656666666666em]

$4_{1} \# 7_{4}$ & $95.56$ & $95.54$ \\
 
$4_{1} \# 7_{4}^m$ & $95.56$ & $95.54$ \\
\addlinespace[0.484656666666666em]

$4_{1} \# 7_{5}$ & $96.37$ & $96.35$ \\
 
$4_{1} \# 7_{5}^m$ & $96.37$ & $96.35$ \\
\addlinespace[0.484656666666666em]

$4_{1} \# 7_{6}$ & $96.35$ & $96.33$ \\
 
$4_{1} \# 7_{6}^m$ & $96.35$ & $96.33$ \\
\addlinespace[0.484656666666666em]

$4_{1} \# 7_{7}$ & $96.76$ & $96.75$ \\
 
$4_{1} \# 7_{7}^m$ & $96.76$ & $96.75$ \\
\addlinespace[0.484656666666666em]

$5_{1} \# 6_{1}$ & $94.62$ & $94.61$ \\
 
$5_{1} \# 6_{1}^m$ & $94.17$ & $94.14$ \\
 
$5_{1}^m \# 6_{1}$ & $94.17$ & $94.14$ \\
 
$5_{1}^m \# 6_{1}^m$ & $94.62$ & $94.61$ \\
\addlinespace[0.484656666666666em]

$5_{1} \# 6_{2}$ & $94.58$ & $94.56$ \\
 
$5_{1} \# 6_{2}^m$ & $95.19$ & $95.17$ \\
 
$5_{1}^m \# 6_{2}$ & $95.19$ & $95.17$ \\
 
$5_{1}^m \# 6_{2}^m$ & $94.58$ & $94.56$ \\
\addlinespace[0.484656666666666em]

$5_{1} \# 6_{3}$ & $95.43$ & $95.41$ \\
 
$5_{1}^m \# 6_{3}$ & $95.43$ & $95.41$ \\
\addlinespace[0.484656666666666em]

$5_{2} \# 6_{1}$ & $95.73$ & $95.71$ \\
 
$5_{2} \# 6_{1}^m$ & $96.36$ & $96.34$ \\

\bottomrule
\end{tabular}
\hfill\hfill\end{center}
\end{table}

\newpage
\begin{table}[h]
 \caption{\label{ByKnotTable 2} Ropelengths of Tight Knots  by Knot Type, Part 3 of 6}

\begin{center}
\hfill\begin{tabular}[t]{l@{\hspace{0.1in}}r@{\hspace{0.1in}}r} \toprule
Knot & $\PRop$ & $\Rop$\\ \midrule
$5_{2}^m \# 6_{1}$ & $96.36$ & $96.34$ \\
 
$5_{2}^m \# 6_{1}^m$ & $95.73$ & $95.71$ \\
\addlinespace[0.523639444444445em]

$5_{2} \# 6_{2}$ & $96.56$ & $96.55$ \\
 
$5_{2} \# 6_{2}^m$ & $95.64$ & $95.62$ \\
 
$5_{2}^m \# 6_{2}$ & $95.64$ & $95.62$ \\
 
$5_{2}^m \# 6_{2}^m$ & $96.56$ & $96.55$ \\
\addlinespace[0.523639444444445em]

$5_{2} \# 6_{3}$ & $96.32$ & $96.30$ \\
 
$5_{2}^m \# 6_{3}$ & $96.32$ & $96.30$ \\
\addlinespace[0.523639444444445em]\cmidrule(r{1.5em}l{1.5em}){1-3}\addlinespace[0.523639444444445em]

$3_{1} \# 3_{1} \# 5_{1}$ & $93.79$ & $93.77$ \\
 
$3_{1} \# 3_{1} \# 5_{1}^m$ & $94.12$ & $94.10$ \\
 
$3_{1} \# 3_{1}^m \# 5_{1}$ & $94.08$ & $94.06$ \\
 
$3_{1} \# 3_{1}^m \# 5_{1}^m$ & $94.08$ & $94.06$ \\
 
$3_{1}^m \# 3_{1}^m \# 5_{1}$ & $94.12$ & $94.10$ \\
 
$3_{1}^m \# 3_{1}^m \# 5_{1}^m$ & $93.79$ & $93.77$ \\
\addlinespace[0.523639444444445em]

$3_{1} \# 3_{1} \# 5_{2}$ & $94.81$ & $94.79$ \\
 
$3_{1} \# 3_{1} \# 5_{2}^m$ & $95.04$ & $95.02$ \\
 
$3_{1} \# 3_{1}^m \# 5_{2}$ & $95.13$ & $95.08$ \\
 
$3_{1} \# 3_{1}^m \# 5_{2}^m$ & $95.13$ & $95.08$ \\
 
$3_{1}^m \# 3_{1}^m \# 5_{2}$ & $95.04$ & $95.02$ \\
 
$3_{1}^m \# 3_{1}^m \# 5_{2}^m$ & $94.81$ & $94.79$ \\
\addlinespace[0.523639444444445em]

$3_{1} \# 4_{1} \# 4_{1}$ & $96.15$ & $96.13$ \\
 
$3_{1}^m \# 4_{1} \# 4_{1}$ & $96.15$ & $96.13$ \\
\addlinespace[0.523639444444445em]\cmidrule(r{0.5em}l{0.5em}){1-3}\addlinespace[0.523639444444445em]

$3_{1} \# 9_{1}$ & $99.88$ & $99.86$ \\
 
$3_{1} \# 9_{1}^m$ & $99.90$ & $99.88$ \\
 
$3_{1}^m \# 9_{1}$ & $99.90$ & $99.88$ \\
 
$3_{1}^m \# 9_{1}^m$ & $99.88$ & $99.86$ \\
\addlinespace[0.523639444444445em]

$3_{1} \# 9_{2}$ & $101.10$ & $101.08$ \\
 
$3_{1} \# 9_{2}^m$ & $101.20$ & $101.17$ \\
 
$3_{1}^m \# 9_{2}$ & $101.20$ & $101.17$ \\
 
$3_{1}^m \# 9_{2}^m$ & $101.10$ & $101.08$ \\

\bottomrule
\end{tabular}
\hfill\begin{tabular}[t]{l@{\hspace{0.1in}}r@{\hspace{0.1in}}r} \toprule
Knot & $\PRop$ & $\Rop$\\ \midrule
$3_{1} \# 9_{3}$ & $101.23$ & $101.21$ \\
 
$3_{1} \# 9_{3}^m$ & $101.27$ & $101.25$ \\
 
$3_{1}^m \# 9_{3}$ & $101.27$ & $101.25$ \\
 
$3_{1}^m \# 9_{3}^m$ & $101.23$ & $101.21$ \\
\addlinespace[0.571125476190476em]

$3_{1} \# 9_{4}$ & $101.32$ & $101.30$ \\
 
$3_{1} \# 9_{4}^m$ & $101.36$ & $101.34$ \\
 
$3_{1}^m \# 9_{4}$ & $101.36$ & $101.34$ \\
 
$3_{1}^m \# 9_{4}^m$ & $101.32$ & $101.30$ \\
\addlinespace[0.571125476190476em]

$3_{1} \# 9_{5}$ & $102.61$ & $102.58$ \\
 
$3_{1} \# 9_{5}^m$ & $102.40$ & $102.37$ \\
 
$3_{1}^m \# 9_{5}$ & $102.40$ & $102.37$ \\
 
$3_{1}^m \# 9_{5}^m$ & $102.61$ & $102.58$ \\
\addlinespace[0.571125476190476em]

$3_{1} \# 9_{6}$ & $102.48$ & $102.45$ \\
 
$3_{1} \# 9_{6}^m$ & $102.64$ & $102.61$ \\
 
$3_{1}^m \# 9_{6}$ & $102.64$ & $102.61$ \\
 
$3_{1}^m \# 9_{6}^m$ & $102.48$ & $102.45$ \\
\addlinespace[0.571125476190476em]

$3_{1} \# 9_{7}$ & $103.28$ & $103.26$ \\
 
$3_{1} \# 9_{7}^m$ & $102.67$ & $102.65$ \\
 
$3_{1}^m \# 9_{7}$ & $102.67$ & $102.65$ \\
 
$3_{1}^m \# 9_{7}^m$ & $103.28$ & $103.26$ \\
\addlinespace[0.571125476190476em]

$3_{1} \# 9_{8}$ & $102.80$ & $102.78$ \\
 
$3_{1} \# 9_{8}^m$ & $103.00$ & $102.98$ \\
 
$3_{1}^m \# 9_{8}$ & $103.00$ & $102.98$ \\
 
$3_{1}^m \# 9_{8}^m$ & $102.80$ & $102.78$ \\
\addlinespace[0.571125476190476em]

$3_{1} \# 9_{9}$ & $102.66$ & $102.63$ \\
 
$3_{1} \# 9_{9}^m$ & $102.78$ & $102.75$ \\
 
$3_{1}^m \# 9_{9}$ & $102.78$ & $102.75$ \\
 
$3_{1}^m \# 9_{9}^m$ & $102.66$ & $102.63$ \\
\addlinespace[0.571125476190476em]

$3_{1} \# 9_{10}$ & $102.82$ & $102.79$ \\
 
$3_{1} \# 9_{10}^m$ & $103.05$ & $103.03$ \\
 
$3_{1}^m \# 9_{10}$ & $103.05$ & $103.03$ \\

\bottomrule
\end{tabular}
\hfill\begin{tabular}[t]{l@{\hspace{0.1in}}r@{\hspace{0.1in}}r} \toprule
Knot & $\PRop$ & $\Rop$\\ \midrule
$3_{1}^m \# 9_{10}^m$ & $102.82$ & $102.79$ \\
\addlinespace[0.499734999999999em]

$3_{1} \# 9_{11}$ & $102.92$ & $102.90$ \\
 
$3_{1} \# 9_{11}^m$ & $102.81$ & $102.78$ \\
 
$3_{1}^m \# 9_{11}$ & $102.81$ & $102.78$ \\
 
$3_{1}^m \# 9_{11}^m$ & $102.92$ & $102.90$ \\
\addlinespace[0.499734999999999em]

$3_{1} \# 9_{12}$ & $103.70$ & $103.67$ \\
 
$3_{1} \# 9_{12}^m$ & $102.93$ & $102.91$ \\
 
$3_{1}^m \# 9_{12}$ & $102.93$ & $102.91$ \\
 
$3_{1}^m \# 9_{12}^m$ & $103.70$ & $103.67$ \\
\addlinespace[0.499734999999999em]

$3_{1} \# 9_{13}$ & $104.61$ & $104.59$ \\
 
$3_{1} \# 9_{13}^m$ & $104.59$ & $104.57$ \\
 
$3_{1}^m \# 9_{13}$ & $104.59$ & $104.57$ \\
 
$3_{1}^m \# 9_{13}^m$ & $104.61$ & $104.59$ \\
\addlinespace[0.499734999999999em]

$3_{1} \# 9_{14}$ & $103.58$ & $103.56$ \\
 
$3_{1} \# 9_{14}^m$ & $103.21$ & $103.19$ \\
 
$3_{1}^m \# 9_{14}$ & $103.21$ & $103.19$ \\
 
$3_{1}^m \# 9_{14}^m$ & $103.58$ & $103.56$ \\
\addlinespace[0.499734999999999em]

$3_{1} \# 9_{15}$ & $103.61$ & $103.58$ \\
 
$3_{1} \# 9_{15}^m$ & $103.56$ & $103.54$ \\
 
$3_{1}^m \# 9_{15}$ & $103.56$ & $103.54$ \\
 
$3_{1}^m \# 9_{15}^m$ & $103.61$ & $103.58$ \\
\addlinespace[0.499734999999999em]

$3_{1} \# 9_{16}$ & $103.32$ & $103.29$ \\
 
$3_{1} \# 9_{16}^m$ & $103.87$ & $103.85$ \\
 
$3_{1}^m \# 9_{16}$ & $103.87$ & $103.85$ \\
 
$3_{1}^m \# 9_{16}^m$ & $103.32$ & $103.29$ \\
\addlinespace[0.499734999999999em]

$3_{1} \# 9_{17}$ & $103.02$ & $103.00$ \\
 
$3_{1} \# 9_{17}^m$ & $103.65$ & $103.62$ \\
 
$3_{1}^m \# 9_{17}$ & $103.65$ & $103.62$ \\
 
$3_{1}^m \# 9_{17}^m$ & $103.02$ & $103.00$ \\
\addlinespace[0.499734999999999em]

$3_{1} \# 9_{18}$ & $103.77$ & $103.75$ \\
 
$3_{1} \# 9_{18}^m$ & $103.83$ & $103.80$ \\

\bottomrule
\end{tabular}
\hfill\hfill\end{center}
\end{table}

\newpage
\begin{table}[h]
 \caption{\label{ByKnotTable 3} Ropelengths of Tight Knots  by Knot Type, Part 4 of 6}

\begin{center}
\hfill\begin{tabular}[t]{l@{\hspace{0.1in}}r@{\hspace{0.1in}}r} \toprule
Knot & $\PRop$ & $\Rop$\\ \midrule
$3_{1}^m \# 9_{18}$ & $103.83$ & $103.80$ \\
 
$3_{1}^m \# 9_{18}^m$ & $103.77$ & $103.75$ \\
\addlinespace[0.499734999999999em]

$3_{1} \# 9_{19}$ & $103.57$ & $103.55$ \\
 
$3_{1} \# 9_{19}^m$ & $104.08$ & $104.06$ \\
 
$3_{1}^m \# 9_{19}$ & $104.08$ & $104.06$ \\
 
$3_{1}^m \# 9_{19}^m$ & $103.57$ & $103.55$ \\
\addlinespace[0.499734999999999em]

$3_{1} \# 9_{20}$ & $103.67$ & $103.65$ \\
 
$3_{1} \# 9_{20}^m$ & $103.42$ & $103.39$ \\
 
$3_{1}^m \# 9_{20}$ & $103.42$ & $103.39$ \\
 
$3_{1}^m \# 9_{20}^m$ & $103.67$ & $103.65$ \\
\addlinespace[0.499734999999999em]

$3_{1} \# 9_{21}$ & $103.83$ & $103.81$ \\
 
$3_{1} \# 9_{21}^m$ & $103.84$ & $103.82$ \\
 
$3_{1}^m \# 9_{21}$ & $103.84$ & $103.82$ \\
 
$3_{1}^m \# 9_{21}^m$ & $103.83$ & $103.81$ \\
\addlinespace[0.499734999999999em]

$3_{1} \# 9_{22}$ & $103.97$ & $103.94$ \\
 
$3_{1} \# 9_{22}^m$ & $104.12$ & $104.09$ \\
 
$3_{1}^m \# 9_{22}$ & $104.12$ & $104.09$ \\
 
$3_{1}^m \# 9_{22}^m$ & $103.97$ & $103.94$ \\
\addlinespace[0.499734999999999em]

$3_{1} \# 9_{23}$ & $104.59$ & $104.56$ \\
 
$3_{1} \# 9_{23}^m$ & $104.20$ & $104.18$ \\
 
$3_{1}^m \# 9_{23}$ & $104.20$ & $104.18$ \\
 
$3_{1}^m \# 9_{23}^m$ & $104.59$ & $104.56$ \\
\addlinespace[0.499734999999999em]

$3_{1} \# 9_{24}$ & $104.50$ & $104.48$ \\
 
$3_{1} \# 9_{24}^m$ & $103.93$ & $103.91$ \\
 
$3_{1}^m \# 9_{24}$ & $103.93$ & $103.91$ \\
 
$3_{1}^m \# 9_{24}^m$ & $104.50$ & $104.48$ \\
\addlinespace[0.499734999999999em]

$3_{1} \# 9_{25}$ & $104.13$ & $104.11$ \\
 
$3_{1} \# 9_{25}^m$ & $103.94$ & $103.92$ \\
 
$3_{1}^m \# 9_{25}$ & $103.94$ & $103.92$ \\
 
$3_{1}^m \# 9_{25}^m$ & $104.13$ & $104.11$ \\
\addlinespace[0.499734999999999em]

$3_{1} \# 9_{26}$ & $104.28$ & $104.25$ \\

\bottomrule
\end{tabular}
\hfill\begin{tabular}[t]{l@{\hspace{0.1in}}r@{\hspace{0.1in}}r} \toprule
Knot & $\PRop$ & $\Rop$\\ \midrule
$3_{1} \# 9_{26}^m$ & $104.21$ & $104.18$ \\
 
$3_{1}^m \# 9_{26}$ & $104.21$ & $104.18$ \\
 
$3_{1}^m \# 9_{26}^m$ & $104.28$ & $104.25$ \\
\addlinespace[0.666312499999999em]

$3_{1} \# 9_{27}$ & $104.98$ & $104.95$ \\
 
$3_{1} \# 9_{27}^m$ & $104.47$ & $104.44$ \\
 
$3_{1}^m \# 9_{27}$ & $104.47$ & $104.44$ \\
 
$3_{1}^m \# 9_{27}^m$ & $104.98$ & $104.95$ \\
\addlinespace[0.666312499999999em]

$3_{1} \# 9_{28}$ & $104.25$ & $104.23$ \\
 
$3_{1} \# 9_{28}^m$ & $104.52$ & $104.50$ \\
 
$3_{1}^m \# 9_{28}$ & $104.52$ & $104.50$ \\
 
$3_{1}^m \# 9_{28}^m$ & $104.25$ & $104.23$ \\
\addlinespace[0.666312499999999em]

$3_{1} \# 9_{29}$ & $104.13$ & $104.11$ \\
 
$3_{1} \# 9_{29}^m$ & $104.97$ & $104.94$ \\
 
$3_{1}^m \# 9_{29}$ & $104.97$ & $104.94$ \\
 
$3_{1}^m \# 9_{29}^m$ & $104.13$ & $104.11$ \\
\addlinespace[0.666312499999999em]

$3_{1} \# 9_{30}$ & $104.39$ & $104.36$ \\
 
$3_{1} \# 9_{30}^m$ & $104.20$ & $104.18$ \\
 
$3_{1}^m \# 9_{30}$ & $104.20$ & $104.18$ \\
 
$3_{1}^m \# 9_{30}^m$ & $104.39$ & $104.36$ \\
\addlinespace[0.666312499999999em]

$3_{1} \# 9_{31}$ & $104.80$ & $104.78$ \\
 
$3_{1} \# 9_{31}^m$ & $105.35$ & $105.33$ \\
 
$3_{1}^m \# 9_{31}$ & $105.35$ & $105.33$ \\
 
$3_{1}^m \# 9_{31}^m$ & $104.80$ & $104.78$ \\
\addlinespace[0.666312499999999em]

$3_{1} \# 9_{32}$ & $104.64$ & $104.62$ \\
 
$3_{1} \# 9_{32}^r$ & $104.64$ & $104.62$ \\
 
$3_{1} \# 9_{32}^m$ & $104.60$ & $104.58$ \\
 
$3_{1} \# 9_{32}^{rm}$ & $104.60$ & $104.58$ \\
 
$3_{1}^m \# 9_{32}$ & $104.60$ & $104.58$ \\
 
$3_{1}^m \# 9_{32}^r$ & $104.60$ & $104.58$ \\
 
$3_{1}^m \# 9_{32}^m$ & $104.64$ & $104.62$ \\
 
$3_{1}^m \# 9_{32}^{rm}$ & $104.64$ & $104.62$ \\

\bottomrule
\end{tabular}
\hfill\begin{tabular}[t]{l@{\hspace{0.1in}}r@{\hspace{0.1in}}r} \toprule
Knot & $\PRop$ & $\Rop$\\ \midrule
$3_{1} \# 9_{33}$ & $105.77$ & $105.74$ \\
 
$3_{1} \# 9_{33}^r$ & $105.77$ & $105.74$ \\
 
$3_{1} \# 9_{33}^m$ & $104.86$ & $104.83$ \\
 
$3_{1} \# 9_{33}^{rm}$ & $104.86$ & $104.83$ \\
 
$3_{1}^m \# 9_{33}$ & $104.86$ & $104.83$ \\
 
$3_{1}^m \# 9_{33}^r$ & $104.86$ & $104.83$ \\
 
$3_{1}^m \# 9_{33}^m$ & $105.77$ & $105.74$ \\
 
$3_{1}^m \# 9_{33}^{rm}$ & $105.77$ & $105.74$ \\
\addlinespace[0.666312499999999em]

$3_{1} \# 9_{34}$ & $105.07$ & $105.05$ \\
 
$3_{1} \# 9_{34}^m$ & $105.62$ & $105.60$ \\
 
$3_{1}^m \# 9_{34}$ & $105.62$ & $105.60$ \\
 
$3_{1}^m \# 9_{34}^m$ & $105.07$ & $105.05$ \\
\addlinespace[0.666312499999999em]

$3_{1} \# 9_{35}$ & $102.44$ & $102.42$ \\
 
$3_{1} \# 9_{35}^m$ & $102.82$ & $102.80$ \\
 
$3_{1}^m \# 9_{35}$ & $102.82$ & $102.80$ \\
 
$3_{1}^m \# 9_{35}^m$ & $102.44$ & $102.42$ \\
\addlinespace[0.666312499999999em]

$3_{1} \# 9_{36}$ & $103.90$ & $103.87$ \\
 
$3_{1} \# 9_{36}^m$ & $103.40$ & $103.38$ \\
 
$3_{1}^m \# 9_{36}$ & $103.40$ & $103.38$ \\
 
$3_{1}^m \# 9_{36}^m$ & $103.90$ & $103.87$ \\
\addlinespace[0.666312499999999em]

$3_{1} \# 9_{37}$ & $103.54$ & $103.52$ \\
 
$3_{1} \# 9_{37}^m$ & $103.75$ & $103.73$ \\
 
$3_{1}^m \# 9_{37}$ & $103.75$ & $103.73$ \\
 
$3_{1}^m \# 9_{37}^m$ & $103.54$ & $103.52$ \\
\addlinespace[0.666312499999999em]

$3_{1} \# 9_{38}$ & $104.63$ & $104.61$ \\
 
$3_{1} \# 9_{38}^m$ & $104.76$ & $104.74$ \\
 
$3_{1}^m \# 9_{38}$ & $104.76$ & $104.74$ \\
 
$3_{1}^m \# 9_{38}^m$ & $104.63$ & $104.61$ \\
\addlinespace[0.666312499999999em]

$3_{1} \# 9_{39}$ & $104.53$ & $104.51$ \\
 
$3_{1} \# 9_{39}^m$ & $104.85$ & $104.83$ \\
 
$3_{1}^m \# 9_{39}$ & $104.85$ & $104.83$ \\

\bottomrule
\end{tabular}
\hfill\hfill\end{center}
\end{table}

\newpage
\begin{table}[h]
 \caption{\label{ByKnotTable 4} Ropelengths of Tight Knots  by Knot Type, Part 5 of 6}

\begin{center}
\hfill\begin{tabular}[t]{l@{\hspace{0.1in}}r@{\hspace{0.1in}}r} \toprule
Knot & $\PRop$ & $\Rop$\\ \midrule
$3_{1}^m \# 9_{39}^m$ & $104.53$ & $104.51$ \\
\addlinespace[0.499734999999999em]

$3_{1} \# 9_{40}$ & $105.91$ & $105.89$ \\
 
$3_{1} \# 9_{40}^m$ & $105.76$ & $105.74$ \\
 
$3_{1}^m \# 9_{40}$ & $105.76$ & $105.74$ \\
 
$3_{1}^m \# 9_{40}^m$ & $105.91$ & $105.89$ \\
\addlinespace[0.499734999999999em]

$3_{1} \# 9_{41}$ & $104.67$ & $104.65$ \\
 
$3_{1} \# 9_{41}^m$ & $104.39$ & $104.37$ \\
 
$3_{1}^m \# 9_{41}$ & $104.39$ & $104.37$ \\
 
$3_{1}^m \# 9_{41}^m$ & $104.67$ & $104.65$ \\
\addlinespace[0.499734999999999em]

$3_{1} \# 9_{42}$ & $92.02$ & $92.00$ \\
 
$3_{1} \# 9_{42}^m$ & $91.79$ & $91.77$ \\
 
$3_{1}^m \# 9_{42}$ & $91.79$ & $91.77$ \\
 
$3_{1}^m \# 9_{42}^m$ & $92.02$ & $92.00$ \\
\addlinespace[0.499734999999999em]

$3_{1} \# 9_{43}$ & $93.70$ & $93.69$ \\
 
$3_{1} \# 9_{43}^m$ & $93.81$ & $93.78$ \\
 
$3_{1}^m \# 9_{43}$ & $93.81$ & $93.78$ \\
 
$3_{1}^m \# 9_{43}^m$ & $93.70$ & $93.69$ \\
\addlinespace[0.499734999999999em]

$3_{1} \# 9_{44}$ & $93.86$ & $93.84$ \\
 
$3_{1} \# 9_{44}^m$ & $94.10$ & $94.08$ \\
 
$3_{1}^m \# 9_{44}$ & $94.10$ & $94.08$ \\
 
$3_{1}^m \# 9_{44}^m$ & $93.86$ & $93.84$ \\
\addlinespace[0.499734999999999em]

$3_{1} \# 9_{45}$ & $97.26$ & $97.24$ \\
 
$3_{1} \# 9_{45}^m$ & $97.71$ & $97.69$ \\
 
$3_{1}^m \# 9_{45}$ & $97.71$ & $97.69$ \\
 
$3_{1}^m \# 9_{45}^m$ & $97.26$ & $97.24$ \\
\addlinespace[0.499734999999999em]

$3_{1} \# 9_{46}$ & $90.93$ & $90.91$ \\
 
$3_{1} \# 9_{46}^m$ & $91.33$ & $91.31$ \\
 
$3_{1}^m \# 9_{46}$ & $91.33$ & $91.31$ \\
 
$3_{1}^m \# 9_{46}^m$ & $90.93$ & $90.91$ \\
\addlinespace[0.499734999999999em]

$3_{1} \# 9_{47}$ & $97.55$ & $97.51$ \\
 
$3_{1} \# 9_{47}^m$ & $97.17$ & $97.15$ \\

\bottomrule
\end{tabular}
\hfill\begin{tabular}[t]{l@{\hspace{0.1in}}r@{\hspace{0.1in}}r} \toprule
Knot & $\PRop$ & $\Rop$\\ \midrule
$3_{1}^m \# 9_{47}$ & $97.17$ & $97.15$ \\
 
$3_{1}^m \# 9_{47}^m$ & $97.55$ & $97.51$ \\
\addlinespace[0.512658589743589em]

$3_{1} \# 9_{48}$ & $97.13$ & $97.11$ \\
 
$3_{1} \# 9_{48}^m$ & $96.50$ & $96.47$ \\
 
$3_{1}^m \# 9_{48}$ & $96.50$ & $96.47$ \\
 
$3_{1}^m \# 9_{48}^m$ & $97.13$ & $97.11$ \\
\addlinespace[0.512658589743589em]

$3_{1} \# 9_{49}$ & $96.29$ & $96.27$ \\
 
$3_{1} \# 9_{49}^m$ & $96.17$ & $96.16$ \\
 
$3_{1}^m \# 9_{49}$ & $96.17$ & $96.16$ \\
 
$3_{1}^m \# 9_{49}^m$ & $96.29$ & $96.27$ \\
\addlinespace[0.512658589743589em]

$4_{1} \# 8_{1}$ & $101.81$ & $101.78$ \\
 
$4_{1} \# 8_{1}^m$ & $101.81$ & $101.78$ \\
\addlinespace[0.512658589743589em]

$4_{1} \# 8_{2}$ & $102.44$ & $102.42$ \\
 
$4_{1} \# 8_{2}^m$ & $102.44$ & $102.42$ \\
\addlinespace[0.512658589743589em]

$4_{1} \# 8_{3}$ & $103.69$ & $103.67$ \\
\addlinespace[0.512658589743589em]

$4_{1} \# 8_{4}$ & $102.78$ & $102.76$ \\
 
$4_{1} \# 8_{4}^m$ & $102.78$ & $102.76$ \\
\addlinespace[0.512658589743589em]

$4_{1} \# 8_{5}$ & $103.09$ & $103.07$ \\
 
$4_{1} \# 8_{5}^m$ & $103.09$ & $103.07$ \\
\addlinespace[0.512658589743589em]

$4_{1} \# 8_{6}$ & $103.47$ & $103.45$ \\
 
$4_{1} \# 8_{6}^m$ & $103.47$ & $103.45$ \\
\addlinespace[0.512658589743589em]

$4_{1} \# 8_{7}$ & $103.55$ & $103.53$ \\
 
$4_{1} \# 8_{7}^m$ & $103.55$ & $103.53$ \\
\addlinespace[0.512658589743589em]

$4_{1} \# 8_{8}$ & $103.50$ & $103.48$ \\
 
$4_{1} \# 8_{8}^m$ & $103.50$ & $103.48$ \\
\addlinespace[0.512658589743589em]

$4_{1} \# 8_{9}$ & $103.94$ & $103.93$ \\
\addlinespace[0.512658589743589em]

$4_{1} \# 8_{10}$ & $104.42$ & $104.40$ \\
 
$4_{1} \# 8_{10}^m$ & $104.42$ & $104.40$ \\
\addlinespace[0.512658589743589em]

$4_{1} \# 8_{11}$ & $103.76$ & $103.74$ \\

\bottomrule
\end{tabular}
\hfill\begin{tabular}[t]{l@{\hspace{0.1in}}r@{\hspace{0.1in}}r} \toprule
Knot & $\PRop$ & $\Rop$\\ \midrule
$4_{1} \# 8_{11}^m$ & $103.76$ & $103.74$ \\
\addlinespace[0.512658589743589em]

$4_{1} \# 8_{12}$ & $103.87$ & $103.84$ \\
\addlinespace[0.512658589743589em]

$4_{1} \# 8_{13}$ & $104.19$ & $104.17$ \\
 
$4_{1} \# 8_{13}^m$ & $104.19$ & $104.17$ \\
\addlinespace[0.512658589743589em]

$4_{1} \# 8_{14}$ & $104.33$ & $104.30$ \\
 
$4_{1} \# 8_{14}^m$ & $104.33$ & $104.30$ \\
\addlinespace[0.512658589743589em]

$4_{1} \# 8_{15}$ & $104.23$ & $104.21$ \\
 
$4_{1} \# 8_{15}^m$ & $104.23$ & $104.21$ \\
\addlinespace[0.512658589743589em]

$4_{1} \# 8_{16}$ & $104.74$ & $104.71$ \\
 
$4_{1} \# 8_{16}^m$ & $104.74$ & $104.71$ \\
\addlinespace[0.512658589743589em]

$4_{1} \# 8_{17}$ & $104.96$ & $104.94$ \\
 
$4_{1} \# 8_{17}^r$ & $104.96$ & $104.94$ \\
\addlinespace[0.512658589743589em]

$4_{1} \# 8_{18}$ & $106.53$ & $106.51$ \\
\addlinespace[0.512658589743589em]

$4_{1} \# 8_{19}$ & $91.64$ & $91.62$ \\
 
$4_{1} \# 8_{19}^m$ & $91.64$ & $91.62$ \\
\addlinespace[0.512658589743589em]

$4_{1} \# 8_{20}$ & $93.64$ & $93.62$ \\
 
$4_{1} \# 8_{20}^m$ & $93.64$ & $93.62$ \\
\addlinespace[0.512658589743589em]

$4_{1} \# 8_{21}$ & $96.20$ & $96.19$ \\
 
$4_{1} \# 8_{21}^m$ & $96.20$ & $96.19$ \\
\addlinespace[0.512658589743589em]

$5_{1} \# 7_{1}$ & $100.38$ & $100.36$ \\
 
$5_{1} \# 7_{1}^m$ & $100.45$ & $100.43$ \\
 
$5_{1}^m \# 7_{1}$ & $100.45$ & $100.43$ \\
 
$5_{1}^m \# 7_{1}^m$ & $100.38$ & $100.36$ \\
\addlinespace[0.512658589743589em]

$5_{1} \# 7_{2}$ & $101.30$ & $101.27$ \\
 
$5_{1} \# 7_{2}^m$ & $101.14$ & $101.12$ \\
 
$5_{1}^m \# 7_{2}$ & $101.14$ & $101.12$ \\
 
$5_{1}^m \# 7_{2}^m$ & $101.30$ & $101.27$ \\
\addlinespace[0.512658589743589em]

$5_{1} \# 7_{3}$ & $101.36$ & $101.34$ \\
 
$5_{1} \# 7_{3}^m$ & $101.32$ & $101.30$ \\

\bottomrule
\end{tabular}
\hfill\hfill\end{center}
\end{table}

\newpage
\begin{table}[h]
 \caption{\label{ByKnotTable 5} Ropelengths of Tight Knots  by Knot Type, Part 6 of 6}

\begin{center}
\hfill\begin{tabular}[t]{l@{\hspace{0.1in}}r@{\hspace{0.1in}}r} \toprule
Knot & $\PRop$ & $\Rop$\\ \midrule
$5_{1}^m \# 7_{3}$ & $101.32$ & $101.30$ \\
 
$5_{1}^m \# 7_{3}^m$ & $101.36$ & $101.34$ \\
\addlinespace[0.499734999999999em]

$5_{1} \# 7_{4}$ & $102.29$ & $102.27$ \\
 
$5_{1} \# 7_{4}^m$ & $101.88$ & $101.86$ \\
 
$5_{1}^m \# 7_{4}$ & $101.88$ & $101.86$ \\
 
$5_{1}^m \# 7_{4}^m$ & $102.29$ & $102.27$ \\
\addlinespace[0.499734999999999em]

$5_{1} \# 7_{5}$ & $102.61$ & $102.59$ \\
 
$5_{1} \# 7_{5}^m$ & $102.79$ & $102.77$ \\
 
$5_{1}^m \# 7_{5}$ & $102.79$ & $102.77$ \\
 
$5_{1}^m \# 7_{5}^m$ & $102.61$ & $102.59$ \\
\addlinespace[0.499734999999999em]

$5_{1} \# 7_{6}$ & $102.97$ & $102.95$ \\
 
$5_{1} \# 7_{6}^m$ & $102.80$ & $102.78$ \\
 
$5_{1}^m \# 7_{6}$ & $102.80$ & $102.78$ \\
 
$5_{1}^m \# 7_{6}^m$ & $102.97$ & $102.95$ \\
\addlinespace[0.499734999999999em]

$5_{1} \# 7_{7}$ & $104.03$ & $104.00$ \\
 
$5_{1} \# 7_{7}^m$ & $103.14$ & $103.12$ \\
 
$5_{1}^m \# 7_{7}$ & $103.14$ & $103.12$ \\
 
$5_{1}^m \# 7_{7}^m$ & $104.03$ & $104.00$ \\
\addlinespace[0.499734999999999em]

$5_{2} \# 7_{1}$ & $101.01$ & $100.99$ \\
 
$5_{2} \# 7_{1}^m$ & $101.19$ & $101.16$ \\
 
$5_{2}^m \# 7_{1}$ & $101.19$ & $101.16$ \\
 
$5_{2}^m \# 7_{1}^m$ & $101.01$ & $100.99$ \\
\addlinespace[0.499734999999999em]

$5_{2} \# 7_{2}$ & $102.55$ & $102.53$ \\
 
$5_{2} \# 7_{2}^m$ & $103.66$ & $103.64$ \\
 
$5_{2}^m \# 7_{2}$ & $103.66$ & $103.64$ \\
 
$5_{2}^m \# 7_{2}^m$ & $102.55$ & $102.53$ \\
\addlinespace[0.499734999999999em]

$5_{2} \# 7_{3}$ & $102.22$ & $102.21$ \\
 
$5_{2} \# 7_{3}^m$ & $102.19$ & $102.17$ \\
 
$5_{2}^m \# 7_{3}$ & $102.19$ & $102.17$ \\
 
$5_{2}^m \# 7_{3}^m$ & $102.22$ & $102.21$ \\
\addlinespace[0.499734999999999em]

$5_{2} \# 7_{4}$ & $102.96$ & $102.94$ \\

\bottomrule
\end{tabular}
\hfill\begin{tabular}[t]{l@{\hspace{0.1in}}r@{\hspace{0.1in}}r} \toprule
Knot & $\PRop$ & $\Rop$\\ \midrule
$5_{2} \# 7_{4}^m$ & $102.76$ & $102.74$ \\
 
$5_{2}^m \# 7_{4}$ & $102.76$ & $102.74$ \\
 
$5_{2}^m \# 7_{4}^m$ & $102.96$ & $102.94$ \\
\addlinespace[0.592357407407407em]

$5_{2} \# 7_{5}$ & $103.39$ & $103.36$ \\
 
$5_{2} \# 7_{5}^m$ & $103.36$ & $103.34$ \\
 
$5_{2}^m \# 7_{5}$ & $103.36$ & $103.34$ \\
 
$5_{2}^m \# 7_{5}^m$ & $103.39$ & $103.36$ \\
\addlinespace[0.592357407407407em]

$5_{2} \# 7_{6}$ & $103.57$ & $103.55$ \\
 
$5_{2} \# 7_{6}^m$ & $103.51$ & $103.49$ \\
 
$5_{2}^m \# 7_{6}$ & $103.51$ & $103.49$ \\
 
$5_{2}^m \# 7_{6}^m$ & $103.57$ & $103.55$ \\
\addlinespace[0.592357407407407em]

$5_{2} \# 7_{7}$ & $105.36$ & $105.34$ \\
 
$5_{2} \# 7_{7}^m$ & $103.93$ & $103.91$ \\
 
$5_{2}^m \# 7_{7}$ & $103.93$ & $103.91$ \\
 
$5_{2}^m \# 7_{7}^m$ & $105.36$ & $105.34$ \\
\addlinespace[0.592357407407407em]

$6_{1} \# 6_{1}$ & $102.44$ & $102.42$ \\
 
$6_{1} \# 6_{1}^m$ & $105.16$ & $105.14$ \\
 
$6_{1}^m \# 6_{1}^m$ & $102.44$ & $102.42$ \\
\addlinespace[0.592357407407407em]

$6_{1} \# 6_{2}$ & $102.53$ & $102.51$ \\
 
$6_{1} \# 6_{2}^m$ & $104.19$ & $104.17$ \\
 
$6_{1}^m \# 6_{2}$ & $104.19$ & $104.17$ \\
 
$6_{1}^m \# 6_{2}^m$ & $102.53$ & $102.51$ \\
\addlinespace[0.592357407407407em]

$6_{1} \# 6_{3}$ & $103.57$ & $103.54$ \\
 
$6_{1}^m \# 6_{3}$ & $103.57$ & $103.54$ \\
\addlinespace[0.592357407407407em]

$6_{2} \# 6_{2}$ & $103.17$ & $103.14$ \\
 
$6_{2} \# 6_{2}^m$ & $103.73$ & $103.71$ \\
 
$6_{2}^m \# 6_{2}^m$ & $103.17$ & $103.14$ \\
\addlinespace[0.592357407407407em]

$6_{2} \# 6_{3}$ & $103.80$ & $103.77$ \\
 
$6_{2}^m \# 6_{3}$ & $103.80$ & $103.77$ \\
\addlinespace[0.592357407407407em]

$6_{3} \# 6_{3}$ & $105.35$ & $105.33$ \\

\bottomrule
\end{tabular}
\hfill\begin{tabular}[t]{l@{\hspace{0.1in}}r@{\hspace{0.1in}}r} \toprule
Knot & $\PRop$ & $\Rop$\\ \midrule
$3_{1} \# 3_{1} \# 6_{1}$ & $102.09$ & $102.06$ \\
 
$3_{1} \# 3_{1} \# 6_{1}^m$ & $103.30$ & $103.28$ \\
 
$3_{1} \# 3_{1}^m \# 6_{1}$ & $102.27$ & $102.24$ \\
 
$3_{1} \# 3_{1}^m \# 6_{1}^m$ & $102.27$ & $102.24$ \\
 
$3_{1}^m \# 3_{1}^m \# 6_{1}$ & $103.30$ & $103.28$ \\
 
$3_{1}^m \# 3_{1}^m \# 6_{1}^m$ & $102.09$ & $102.06$ \\
\addlinespace[0.45em]

$3_{1} \# 3_{1} \# 6_{2}$ & $102.61$ & $102.58$ \\
 
$3_{1} \# 3_{1} \# 6_{2}^m$ & $103.50$ & $103.47$ \\
 
$3_{1} \# 3_{1}^m \# 6_{2}$ & $103.48$ & $103.46$ \\
 
$3_{1} \# 3_{1}^m \# 6_{2}^m$ & $103.48$ & $103.46$ \\
 
$3_{1}^m \# 3_{1}^m \# 6_{2}$ & $103.50$ & $103.47$ \\
 
$3_{1}^m \# 3_{1}^m \# 6_{2}^m$ & $102.61$ & $102.58$ \\
\addlinespace[0.45em]

$3_{1} \# 3_{1} \# 6_{3}$ & $103.69$ & $103.67$ \\
 
$3_{1} \# 3_{1}^m \# 6_{3}$ & $104.01$ & $103.98$ \\
 
$3_{1}^m \# 3_{1}^m \# 6_{3}$ & $103.69$ & $103.67$ \\
\addlinespace[0.45em]

$3_{1} \# 4_{1} \# 5_{1}$ & $102.42$ & $102.40$ \\
 
$3_{1} \# 4_{1} \# 5_{1}^m$ & $102.34$ & $102.31$ \\
 
$3_{1}^m \# 4_{1} \# 5_{1}$ & $102.34$ & $102.31$ \\
 
$3_{1}^m \# 4_{1} \# 5_{1}^m$ & $102.42$ & $102.40$ \\
\addlinespace[0.45em]

$3_{1} \# 4_{1} \# 5_{2}$ & $106.20$ & $106.18$ \\
 
$3_{1} \# 4_{1} \# 5_{2}^m$ & $103.21$ & $103.19$ \\
 
$3_{1}^m \# 4_{1} \# 5_{2}$ & $103.21$ & $103.19$ \\
 
$3_{1}^m \# 4_{1} \# 5_{2}^m$ & $106.20$ & $106.18$ \\
\addlinespace[0.45em]

$4_{1} \# 4_{1} \# 4_{1}$ & $104.16$ & $104.14$ \\
\addlinespace[0.45em]\cmidrule(r{1.5em}l{1.5em}){1-3}\addlinespace[0.45em]

$3_{1} \# 3_{1} \# 3_{1} \# 3_{1}$ & $101.60$ & $101.58$ \\
 
$3_{1} \# 3_{1} \# 3_{1} \# 3_{1}^m$ & $102.80$ & $102.78$ \\
 
$3_{1} \# 3_{1} \# 3_{1}^m \# 3_{1}^m$ & $102.24$ & $102.22$ \\
 
$3_{1} \# 3_{1}^m \# 3_{1}^m \# 3_{1}^m$ & $102.80$ & $102.78$ \\
 
$3_{1}^m \# 3_{1}^m \# 3_{1}^m \# 3_{1}^m$ & $101.60$ & $101.58$ \\

\bottomrule
\end{tabular}
\hfill\hfill\end{center}
\end{table}

\begin{table}[h]
 \caption{\label{ResidualTable 0} Residuals of Tight Knots  by Knot Type, Part 1 of 6}

\begin{center}
\hfill\begin{tabular}[t]{l@{\hspace{0.1in}}r} \toprule
Knot & \hspace{-5pt}Residual\\ \midrule
$3_{1} \# 3_{1}$ & $0.98 \times 10^{-3}$ \\
 
$3_{1} \# 3_{1}^m$ & $2.70 \times 10^{-3}$ \\
 
$3_{1}^m \# 3_{1}^m$ & $0.98 \times 10^{-3}$ \\
\addlinespace[0.589183333333334em]\cmidrule(r{0.5em}l{0.5em}){1-2}\addlinespace[0.589183333333334em]

$3_{1} \# 4_{1}$ & $0.75 \times 10^{-3}$ \\
 
$3_{1}^m \# 4_{1}$ & $0.75 \times 10^{-3}$ \\
\addlinespace[0.589183333333334em]\cmidrule(r{0.5em}l{0.5em}){1-2}\addlinespace[0.589183333333334em]

$3_{1} \# 5_{1}$ & $0.98 \times 10^{-3}$ \\
 
$3_{1} \# 5_{1}^m$ & $2.63 \times 10^{-3}$ \\
 
$3_{1}^m \# 5_{1}$ & $2.63 \times 10^{-3}$ \\
 
$3_{1}^m \# 5_{1}^m$ & $0.98 \times 10^{-3}$ \\
\addlinespace[0.589183333333334em]

$3_{1} \# 5_{2}$ & $0.98 \times 10^{-3}$ \\
 
$3_{1} \# 5_{2}^m$ & $0.75 \times 10^{-3}$ \\
 
$3_{1}^m \# 5_{2}$ & $0.75 \times 10^{-3}$ \\
 
$3_{1}^m \# 5_{2}^m$ & $0.98 \times 10^{-3}$ \\
\addlinespace[0.589183333333334em]

$4_{1} \# 4_{1}$ & $0.75 \times 10^{-3}$ \\
\addlinespace[0.589183333333334em]\cmidrule(r{0.5em}l{0.5em}){1-2}\addlinespace[0.589183333333334em]

$3_{1} \# 6_{1}$ & $4.38 \times 10^{-3}$ \\
 
$3_{1} \# 6_{1}^m$ & $1.12 \times 10^{-3}$ \\
 
$3_{1}^m \# 6_{1}$ & $1.12 \times 10^{-3}$ \\
 
$3_{1}^m \# 6_{1}^m$ & $4.38 \times 10^{-3}$ \\
\addlinespace[0.589183333333334em]

$3_{1} \# 6_{2}$ & $18.72 \times 10^{-3}$ \\
 
$3_{1} \# 6_{2}^m$ & $2.19 \times 10^{-3}$ \\
 
$3_{1}^m \# 6_{2}$ & $2.19 \times 10^{-3}$ \\
 
$3_{1}^m \# 6_{2}^m$ & $18.72 \times 10^{-3}$ \\
\addlinespace[0.589183333333334em]

$3_{1} \# 6_{3}$ & $0.75 \times 10^{-3}$ \\
 
$3_{1}^m \# 6_{3}$ & $0.75 \times 10^{-3}$ \\
\addlinespace[0.589183333333334em]

$4_{1} \# 5_{1}$ & $0.75 \times 10^{-3}$ \\
 
$4_{1} \# 5_{1}^m$ & $0.75 \times 10^{-3}$ \\
\addlinespace[0.589183333333334em]

$4_{1} \# 5_{2}$ & $0.75 \times 10^{-3}$ \\
 
$4_{1} \# 5_{2}^m$ & $0.75 \times 10^{-3}$ \\

\bottomrule
\end{tabular}
\hfill\begin{tabular}[t]{l@{\hspace{0.1in}}r} \toprule
Knot & \hspace{-5pt}Residual\\ \midrule
$3_{1} \# 3_{1} \# 3_{1}$ & $0.75 \times 10^{-3}$ \\
 
$3_{1} \# 3_{1} \# 3_{1}^m$ & $126.45 \times 10^{-3}$ \\
 
$3_{1} \# 3_{1}^m \# 3_{1}^m$ & $126.45 \times 10^{-3}$ \\
 
$3_{1}^m \# 3_{1}^m \# 3_{1}^m$ & $0.75 \times 10^{-3}$ \\
\addlinespace[0.461081041666666em]\cmidrule(r{0.5em}l{0.5em}){1-2}\addlinespace[0.461081041666666em]

$3_{1} \# 7_{1}$ & $0.98 \times 10^{-3}$ \\
 
$3_{1} \# 7_{1}^m$ & $2.29 \times 10^{-3}$ \\
 
$3_{1}^m \# 7_{1}$ & $2.29 \times 10^{-3}$ \\
 
$3_{1}^m \# 7_{1}^m$ & $0.98 \times 10^{-3}$ \\
\addlinespace[0.461081041666666em]

$3_{1} \# 7_{2}$ & $0.75 \times 10^{-3}$ \\
 
$3_{1} \# 7_{2}^m$ & $1.02 \times 10^{-3}$ \\
 
$3_{1}^m \# 7_{2}$ & $1.02 \times 10^{-3}$ \\
 
$3_{1}^m \# 7_{2}^m$ & $0.75 \times 10^{-3}$ \\
\addlinespace[0.461081041666666em]

$3_{1} \# 7_{3}$ & $0.98 \times 10^{-3}$ \\
 
$3_{1} \# 7_{3}^m$ & $0.75 \times 10^{-3}$ \\
 
$3_{1}^m \# 7_{3}$ & $0.75 \times 10^{-3}$ \\
 
$3_{1}^m \# 7_{3}^m$ & $0.98 \times 10^{-3}$ \\
\addlinespace[0.461081041666666em]

$3_{1} \# 7_{4}$ & $8.47 \times 10^{-3}$ \\
 
$3_{1} \# 7_{4}^m$ & $7.50 \times 10^{-3}$ \\
 
$3_{1}^m \# 7_{4}$ & $7.50 \times 10^{-3}$ \\
 
$3_{1}^m \# 7_{4}^m$ & $8.47 \times 10^{-3}$ \\
\addlinespace[0.461081041666666em]

$3_{1} \# 7_{5}$ & $0.75 \times 10^{-3}$ \\
 
$3_{1} \# 7_{5}^m$ & $0.75 \times 10^{-3}$ \\
 
$3_{1}^m \# 7_{5}$ & $0.75 \times 10^{-3}$ \\
 
$3_{1}^m \# 7_{5}^m$ & $0.75 \times 10^{-3}$ \\
\addlinespace[0.461081041666666em]

$3_{1} \# 7_{6}$ & $2.89 \times 10^{-3}$ \\
 
$3_{1} \# 7_{6}^m$ & $0.75 \times 10^{-3}$ \\
 
$3_{1}^m \# 7_{6}$ & $0.75 \times 10^{-3}$ \\
 
$3_{1}^m \# 7_{6}^m$ & $2.89 \times 10^{-3}$ \\
\addlinespace[0.461081041666666em]

$3_{1} \# 7_{7}$ & $3.65 \times 10^{-3}$ \\
 
$3_{1} \# 7_{7}^m$ & $2.71 \times 10^{-3}$ \\
 
$3_{1}^m \# 7_{7}$ & $2.71 \times 10^{-3}$ \\

\bottomrule
\end{tabular}
\hfill\begin{tabular}[t]{l@{\hspace{0.1in}}r} \toprule
Knot & \hspace{-5pt}Residual\\ \midrule
$3_{1}^m \# 7_{7}^m$ & $3.65 \times 10^{-3}$ \\
\addlinespace[0.503841527777778em]

$4_{1} \# 6_{1}$ & $0.75 \times 10^{-3}$ \\
 
$4_{1} \# 6_{1}^m$ & $0.75 \times 10^{-3}$ \\
\addlinespace[0.503841527777778em]

$4_{1} \# 6_{2}$ & $0.75 \times 10^{-3}$ \\
 
$4_{1} \# 6_{2}^m$ & $0.75 \times 10^{-3}$ \\
\addlinespace[0.503841527777778em]

$4_{1} \# 6_{3}$ & $0.75 \times 10^{-3}$ \\
\addlinespace[0.503841527777778em]

$5_{1} \# 5_{1}$ & $0.98 \times 10^{-3}$ \\
 
$5_{1} \# 5_{1}^m$ & $0.98 \times 10^{-3}$ \\
 
$5_{1}^m \# 5_{1}^m$ & $0.98 \times 10^{-3}$ \\
\addlinespace[0.503841527777778em]

$5_{1} \# 5_{2}$ & $0.75 \times 10^{-3}$ \\
 
$5_{1} \# 5_{2}^m$ & $0.75 \times 10^{-3}$ \\
 
$5_{1}^m \# 5_{2}$ & $0.75 \times 10^{-3}$ \\
 
$5_{1}^m \# 5_{2}^m$ & $0.75 \times 10^{-3}$ \\
\addlinespace[0.503841527777778em]

$5_{2} \# 5_{2}$ & $0.75 \times 10^{-3}$ \\
 
$5_{2} \# 5_{2}^m$ & $0.75 \times 10^{-3}$ \\
 
$5_{2}^m \# 5_{2}^m$ & $0.75 \times 10^{-3}$ \\
\addlinespace[0.503841527777778em]\cmidrule(r{1.5em}l{1.5em}){1-2}\addlinespace[0.503841527777778em]

$3_{1} \# 3_{1} \# 4_{1}$ & $0.75 \times 10^{-3}$ \\
 
$3_{1} \# 3_{1}^m \# 4_{1}$ & $7.55 \times 10^{-3}$ \\
 
$3_{1}^m \# 3_{1}^m \# 4_{1}$ & $0.75 \times 10^{-3}$ \\
\addlinespace[0.503841527777778em]\cmidrule(r{0.5em}l{0.5em}){1-2}\addlinespace[0.503841527777778em]

$3_{1} \# 8_{1}$ & $7.49 \times 10^{-3}$ \\
 
$3_{1} \# 8_{1}^m$ & $2.74 \times 10^{-3}$ \\
 
$3_{1}^m \# 8_{1}$ & $2.74 \times 10^{-3}$ \\
 
$3_{1}^m \# 8_{1}^m$ & $7.49 \times 10^{-3}$ \\
\addlinespace[0.503841527777778em]

$3_{1} \# 8_{2}$ & $9.62 \times 10^{-3}$ \\
 
$3_{1} \# 8_{2}^m$ & $26.84 \times 10^{-3}$ \\
 
$3_{1}^m \# 8_{2}$ & $26.84 \times 10^{-3}$ \\
 
$3_{1}^m \# 8_{2}^m$ & $9.62 \times 10^{-3}$ \\
\addlinespace[0.503841527777778em]

$3_{1} \# 8_{3}$ & $7.50 \times 10^{-3}$ \\
 
$3_{1}^m \# 8_{3}$ & $7.50 \times 10^{-3}$ \\

\bottomrule
\end{tabular}
\hfill\hfill\end{center}
\end{table}

\newpage
\begin{table}[h]
 \caption{\label{ResidualTable 1} Residuals of Tight Knots  by Knot Type, Part 2 of 6}

\begin{center}
\hfill\begin{tabular}[t]{l@{\hspace{0.1in}}r} \toprule
Knot & \hspace{-5pt}Residual\\ \midrule
$3_{1} \# 8_{4}$ & $0.75 \times 10^{-3}$ \\
 
$3_{1} \# 8_{4}^m$ & $2.70 \times 10^{-3}$ \\
 
$3_{1}^m \# 8_{4}$ & $2.70 \times 10^{-3}$ \\
 
$3_{1}^m \# 8_{4}^m$ & $0.75 \times 10^{-3}$ \\
\addlinespace[0.499734999999999em]

$3_{1} \# 8_{5}$ & $0.75 \times 10^{-3}$ \\
 
$3_{1} \# 8_{5}^m$ & $7.50 \times 10^{-3}$ \\
 
$3_{1}^m \# 8_{5}$ & $7.50 \times 10^{-3}$ \\
 
$3_{1}^m \# 8_{5}^m$ & $0.75 \times 10^{-3}$ \\
\addlinespace[0.499734999999999em]

$3_{1} \# 8_{6}$ & $0.75 \times 10^{-3}$ \\
 
$3_{1} \# 8_{6}^m$ & $0.75 \times 10^{-3}$ \\
 
$3_{1}^m \# 8_{6}$ & $0.75 \times 10^{-3}$ \\
 
$3_{1}^m \# 8_{6}^m$ & $0.75 \times 10^{-3}$ \\
\addlinespace[0.499734999999999em]

$3_{1} \# 8_{7}$ & $0.75 \times 10^{-3}$ \\
 
$3_{1} \# 8_{7}^m$ & $0.98 \times 10^{-3}$ \\
 
$3_{1}^m \# 8_{7}$ & $0.98 \times 10^{-3}$ \\
 
$3_{1}^m \# 8_{7}^m$ & $0.75 \times 10^{-3}$ \\
\addlinespace[0.499734999999999em]

$3_{1} \# 8_{8}$ & $7.46 \times 10^{-3}$ \\
 
$3_{1} \# 8_{8}^m$ & $0.75 \times 10^{-3}$ \\
 
$3_{1}^m \# 8_{8}$ & $0.75 \times 10^{-3}$ \\
 
$3_{1}^m \# 8_{8}^m$ & $7.46 \times 10^{-3}$ \\
\addlinespace[0.499734999999999em]

$3_{1} \# 8_{9}$ & $7.27 \times 10^{-3}$ \\
 
$3_{1}^m \# 8_{9}$ & $7.27 \times 10^{-3}$ \\
\addlinespace[0.499734999999999em]

$3_{1} \# 8_{10}$ & $3.58 \times 10^{-3}$ \\
 
$3_{1} \# 8_{10}^m$ & $7.50 \times 10^{-3}$ \\
 
$3_{1}^m \# 8_{10}$ & $7.50 \times 10^{-3}$ \\
 
$3_{1}^m \# 8_{10}^m$ & $3.58 \times 10^{-3}$ \\
\addlinespace[0.499734999999999em]

$3_{1} \# 8_{11}$ & $0.75 \times 10^{-3}$ \\
 
$3_{1} \# 8_{11}^m$ & $7.50 \times 10^{-3}$ \\
 
$3_{1}^m \# 8_{11}$ & $7.50 \times 10^{-3}$ \\
 
$3_{1}^m \# 8_{11}^m$ & $0.75 \times 10^{-3}$ \\
\addlinespace[0.499734999999999em]

$3_{1} \# 8_{12}$ & $4.36 \times 10^{-3}$ \\

\bottomrule
\end{tabular}
\hfill\begin{tabular}[t]{l@{\hspace{0.1in}}r} \toprule
Knot & \hspace{-5pt}Residual\\ \midrule
$3_{1}^m \# 8_{12}$ & $4.36 \times 10^{-3}$ \\
\addlinespace[0.499734999999999em]

$3_{1} \# 8_{13}$ & $7.50 \times 10^{-3}$ \\
 
$3_{1} \# 8_{13}^m$ & $0.75 \times 10^{-3}$ \\
 
$3_{1}^m \# 8_{13}$ & $0.75 \times 10^{-3}$ \\
 
$3_{1}^m \# 8_{13}^m$ & $7.50 \times 10^{-3}$ \\
\addlinespace[0.499734999999999em]

$3_{1} \# 8_{14}$ & $7.50 \times 10^{-3}$ \\
 
$3_{1} \# 8_{14}^m$ & $0.75 \times 10^{-3}$ \\
 
$3_{1}^m \# 8_{14}$ & $0.75 \times 10^{-3}$ \\
 
$3_{1}^m \# 8_{14}^m$ & $7.50 \times 10^{-3}$ \\
\addlinespace[0.499734999999999em]

$3_{1} \# 8_{15}$ & $2.64 \times 10^{-3}$ \\
 
$3_{1} \# 8_{15}^m$ & $0.75 \times 10^{-3}$ \\
 
$3_{1}^m \# 8_{15}$ & $0.75 \times 10^{-3}$ \\
 
$3_{1}^m \# 8_{15}^m$ & $2.64 \times 10^{-3}$ \\
\addlinespace[0.499734999999999em]

$3_{1} \# 8_{16}$ & $0.75 \times 10^{-3}$ \\
 
$3_{1} \# 8_{16}^m$ & $0.75 \times 10^{-3}$ \\
 
$3_{1}^m \# 8_{16}$ & $0.75 \times 10^{-3}$ \\
 
$3_{1}^m \# 8_{16}^m$ & $0.75 \times 10^{-3}$ \\
\addlinespace[0.499734999999999em]

$3_{1} \# 8_{17}$ & $0.75 \times 10^{-3}$ \\
 
$3_{1} \# 8_{17}^r$ & $0.75 \times 10^{-3}$ \\
 
$3_{1}^m \# 8_{17}$ & $0.75 \times 10^{-3}$ \\
 
$3_{1}^m \# 8_{17}^r$ & $0.75 \times 10^{-3}$ \\
\addlinespace[0.499734999999999em]

$3_{1} \# 8_{18}$ & $0.75 \times 10^{-3}$ \\
 
$3_{1}^m \# 8_{18}$ & $0.75 \times 10^{-3}$ \\
\addlinespace[0.499734999999999em]

$3_{1} \# 8_{19}$ & $5.59 \times 10^{-3}$ \\
 
$3_{1} \# 8_{19}^m$ & $1.78 \times 10^{-3}$ \\
 
$3_{1}^m \# 8_{19}$ & $1.78 \times 10^{-3}$ \\
 
$3_{1}^m \# 8_{19}^m$ & $5.59 \times 10^{-3}$ \\
\addlinespace[0.499734999999999em]

$3_{1} \# 8_{20}$ & $0.98 \times 10^{-3}$ \\
 
$3_{1} \# 8_{20}^m$ & $2.06 \times 10^{-3}$ \\
 
$3_{1}^m \# 8_{20}$ & $2.06 \times 10^{-3}$ \\
 
$3_{1}^m \# 8_{20}^m$ & $0.98 \times 10^{-3}$ \\

\bottomrule
\end{tabular}
\hfill\begin{tabular}[t]{l@{\hspace{0.1in}}r} \toprule
Knot & \hspace{-5pt}Residual\\ \midrule
$3_{1} \# 8_{21}$ & $0.75 \times 10^{-3}$ \\
 
$3_{1} \# 8_{21}^m$ & $0.75 \times 10^{-3}$ \\
 
$3_{1}^m \# 8_{21}$ & $0.75 \times 10^{-3}$ \\
 
$3_{1}^m \# 8_{21}^m$ & $0.75 \times 10^{-3}$ \\
\addlinespace[0.484656666666666em]

$4_{1} \# 7_{1}$ & $0.75 \times 10^{-3}$ \\
 
$4_{1} \# 7_{1}^m$ & $0.75 \times 10^{-3}$ \\
\addlinespace[0.484656666666666em]

$4_{1} \# 7_{2}$ & $0.75 \times 10^{-3}$ \\
 
$4_{1} \# 7_{2}^m$ & $0.75 \times 10^{-3}$ \\
\addlinespace[0.484656666666666em]

$4_{1} \# 7_{3}$ & $0.75 \times 10^{-3}$ \\
 
$4_{1} \# 7_{3}^m$ & $0.75 \times 10^{-3}$ \\
\addlinespace[0.484656666666666em]

$4_{1} \# 7_{4}$ & $1.38 \times 10^{-3}$ \\
 
$4_{1} \# 7_{4}^m$ & $1.38 \times 10^{-3}$ \\
\addlinespace[0.484656666666666em]

$4_{1} \# 7_{5}$ & $2.31 \times 10^{-3}$ \\
 
$4_{1} \# 7_{5}^m$ & $2.31 \times 10^{-3}$ \\
\addlinespace[0.484656666666666em]

$4_{1} \# 7_{6}$ & $2.56 \times 10^{-3}$ \\
 
$4_{1} \# 7_{6}^m$ & $2.56 \times 10^{-3}$ \\
\addlinespace[0.484656666666666em]

$4_{1} \# 7_{7}$ & $0.75 \times 10^{-3}$ \\
 
$4_{1} \# 7_{7}^m$ & $0.75 \times 10^{-3}$ \\
\addlinespace[0.484656666666666em]

$5_{1} \# 6_{1}$ & $0.75 \times 10^{-3}$ \\
 
$5_{1} \# 6_{1}^m$ & $5.50 \times 10^{-3}$ \\
 
$5_{1}^m \# 6_{1}$ & $5.50 \times 10^{-3}$ \\
 
$5_{1}^m \# 6_{1}^m$ & $0.75 \times 10^{-3}$ \\
\addlinespace[0.484656666666666em]

$5_{1} \# 6_{2}$ & $3.14 \times 10^{-3}$ \\
 
$5_{1} \# 6_{2}^m$ & $3.30 \times 10^{-3}$ \\
 
$5_{1}^m \# 6_{2}$ & $3.30 \times 10^{-3}$ \\
 
$5_{1}^m \# 6_{2}^m$ & $3.14 \times 10^{-3}$ \\
\addlinespace[0.484656666666666em]

$5_{1} \# 6_{3}$ & $0.75 \times 10^{-3}$ \\
 
$5_{1}^m \# 6_{3}$ & $0.75 \times 10^{-3}$ \\
\addlinespace[0.484656666666666em]

$5_{2} \# 6_{1}$ & $7.50 \times 10^{-3}$ \\
 
$5_{2} \# 6_{1}^m$ & $7.50 \times 10^{-3}$ \\

\bottomrule
\end{tabular}
\hfill\hfill\end{center}
\end{table}

\newpage
\begin{table}[h]
 \caption{\label{ResidualTable 2} Residuals of Tight Knots  by Knot Type, Part 3 of 6}

\begin{center}
\hfill\begin{tabular}[t]{l@{\hspace{0.1in}}r} \toprule
Knot & \hspace{-5pt}Residual\\ \midrule
$5_{2}^m \# 6_{1}$ & $7.50 \times 10^{-3}$ \\
 
$5_{2}^m \# 6_{1}^m$ & $7.50 \times 10^{-3}$ \\
\addlinespace[0.523639444444445em]

$5_{2} \# 6_{2}$ & $0.75 \times 10^{-3}$ \\
 
$5_{2} \# 6_{2}^m$ & $0.75 \times 10^{-3}$ \\
 
$5_{2}^m \# 6_{2}$ & $0.75 \times 10^{-3}$ \\
 
$5_{2}^m \# 6_{2}^m$ & $0.75 \times 10^{-3}$ \\
\addlinespace[0.523639444444445em]

$5_{2} \# 6_{3}$ & $0.75 \times 10^{-3}$ \\
\bottomrule
\end{tabular}
\hfill\begin{tabular}[t]{l@{\hspace{0.1in}}r} \toprule
Knot & \hspace{-5pt}Residual\\ \midrule
$3_{1} \# 9_{3}$ & $0.75 \times 10^{-3}$ \\
 
$3_{1} \# 9_{3}^m$ & $2.95 \times 10^{-3}$ \\
 
$3_{1}^m \# 9_{3}$ & $2.95 \times 10^{-3}$ \\
 
$3_{1}^m \# 9_{3}^m$ & $0.75 \times 10^{-3}$ \\
\addlinespace[0.571125476190476em]

$3_{1} \# 9_{4}$ & $15.25 \times 10^{-3}$ \\
 
$3_{1} \# 9_{4}^m$ & $0.75 \times 10^{-3}$ \\
 
$3_{1}^m \# 9_{4}$ & $0.75 \times 10^{-3}$ \\
 
$3_{1}^m \# 9_{4}^m$ & $15.25 \times 10^{-3}$ \\
\addlinespace[0.571125476190476em]

$3_{1} \# 9_{5}$ & $0.75 \times 10^{-3}$ \\
 
$3_{1} \# 9_{5}^m$ & $7.82 \times 10^{-3}$ \\
 
$3_{1}^m \# 9_{5}$ & $7.82 \times 10^{-3}$ \\
 
$3_{1}^m \# 9_{5}^m$ & $0.75 \times 10^{-3}$ \\
\addlinespace[0.571125476190476em]

$3_{1} \# 9_{6}$ & $3.43 \times 10^{-3}$ \\
 
$3_{1} \# 9_{6}^m$ & $0.75 \times 10^{-3}$ \\
 
$3_{1}^m \# 9_{6}$ & $0.75 \times 10^{-3}$ \\
 
$3_{1}^m \# 9_{6}^m$ & $3.43 \times 10^{-3}$ \\
\addlinespace[0.571125476190476em]

$3_{1} \# 9_{7}$ & $75.58 \times 10^{-3}$ \\
 
$3_{1} \# 9_{7}^m$ & $0.75 \times 10^{-3}$ \\
 
$3_{1}^m \# 9_{7}$ & $0.75 \times 10^{-3}$ \\
 
$3_{1}^m \# 9_{7}^m$ & $75.58 \times 10^{-3}$ \\
\addlinespace[0.571125476190476em]

$3_{1} \# 9_{8}$ & $7.50 \times 10^{-3}$ \\
 
$3_{1} \# 9_{8}^m$ & $2.44 \times 10^{-3}$ \\
 
$3_{1}^m \# 9_{8}$ & $2.44 \times 10^{-3}$ \\
 
$3_{1}^m \# 9_{8}^m$ & $7.50 \times 10^{-3}$ \\
\addlinespace[0.571125476190476em]

$3_{1} \# 9_{9}$ & $0.98 \times 10^{-3}$ \\
 
$3_{1} \# 9_{9}^m$ & $7.50 \times 10^{-3}$ \\
 
$3_{1}^m \# 9_{9}$ & $7.50 \times 10^{-3}$ \\
 
$3_{1}^m \# 9_{9}^m$ & $0.98 \times 10^{-3}$ \\
\addlinespace[0.571125476190476em]

$3_{1} \# 9_{10}$ & $4.10 \times 10^{-3}$ \\
 
$3_{1} \# 9_{10}^m$ & $0.75 \times 10^{-3}$ \\
 
$3_{1}^m \# 9_{10}$ & $0.75 \times 10^{-3}$ \\

\bottomrule
\end{tabular}
\hfill\begin{tabular}[t]{l@{\hspace{0.1in}}r} \toprule
Knot & \hspace{-5pt}Residual\\ \midrule
$3_{1}^m \# 9_{10}^m$ & $4.10 \times 10^{-3}$ \\
\addlinespace[0.499734999999999em]

$3_{1} \# 9_{11}$ & $7.49 \times 10^{-3}$ \\
 
$3_{1} \# 9_{11}^m$ & $0.75 \times 10^{-3}$ \\
 
$3_{1}^m \# 9_{11}$ & $0.75 \times 10^{-3}$ \\
 
$3_{1}^m \# 9_{11}^m$ & $7.49 \times 10^{-3}$ \\
\addlinespace[0.499734999999999em]

$3_{1} \# 9_{12}$ & $2.59 \times 10^{-3}$ \\
 
$3_{1} \# 9_{12}^m$ & $0.75 \times 10^{-3}$ \\
 
$3_{1}^m \# 9_{12}$ & $0.75 \times 10^{-3}$ \\
 
$3_{1}^m \# 9_{12}^m$ & $2.59 \times 10^{-3}$ \\
\addlinespace[0.499734999999999em]

$3_{1} \# 9_{13}$ & $2.11 \times 10^{-3}$ \\
 
$3_{1} \# 9_{13}^m$ & $1.30 \times 10^{-3}$ \\
 
$3_{1}^m \# 9_{13}$ & $1.30 \times 10^{-3}$ \\
 
$3_{1}^m \# 9_{13}^m$ & $2.11 \times 10^{-3}$ \\
\addlinespace[0.499734999999999em]

$3_{1} \# 9_{14}$ & $7.49 \times 10^{-3}$ \\
 
$3_{1} \# 9_{14}^m$ & $0.75 \times 10^{-3}$ \\
 
$3_{1}^m \# 9_{14}$ & $0.75 \times 10^{-3}$ \\
 
$3_{1}^m \# 9_{14}^m$ & $7.49 \times 10^{-3}$ \\
\addlinespace[0.499734999999999em]

$3_{1} \# 9_{15}$ & $0.75 \times 10^{-3}$ \\
 
$3_{1} \# 9_{15}^m$ & $2.57 \times 10^{-3}$ \\
 
$3_{1}^m \# 9_{15}$ & $2.57 \times 10^{-3}$ \\
 
$3_{1}^m \# 9_{15}^m$ & $0.75 \times 10^{-3}$ \\
\addlinespace[0.499734999999999em]

$3_{1} \# 9_{16}$ & $6.15 \times 10^{-3}$ \\
 
$3_{1} \# 9_{16}^m$ & $7.50 \times 10^{-3}$ \\
 
$3_{1}^m \# 9_{16}$ & $7.50 \times 10^{-3}$ \\
 
$3_{1}^m \# 9_{16}^m$ & $6.15 \times 10^{-3}$ \\
\addlinespace[0.499734999999999em]

$3_{1} \# 9_{17}$ & $2.17 \times 10^{-3}$ \\
 
$3_{1} \# 9_{17}^m$ & $7.50 \times 10^{-3}$ \\
 
$3_{1}^m \# 9_{17}$ & $7.50 \times 10^{-3}$ \\
 
$3_{1}^m \# 9_{17}^m$ & $2.17 \times 10^{-3}$ \\
\addlinespace[0.499734999999999em]

$3_{1} \# 9_{18}$ & $7.50 \times 10^{-3}$ \\
 
$3_{1} \# 9_{18}^m$ & $2.19 \times 10^{-3}$ \\

\bottomrule
\end{tabular}
\hfill\hfill\end{center}
\end{table}

\newpage
\begin{table}[h]
 \caption{\label{ResidualTable 3} Residuals of Tight Knots  by Knot Type, Part 4 of 6}

\begin{center}
\hfill%
\hfill\hfill\end{center}
\end{table}

\end{document}